\documentclass[12pt,a4paper]{article}
\usepackage[cp1251]{inputenc} 
\usepackage[russian]{babel}
\usepackage{amsfonts}

\newcommand{\di}{\displaystyle}
\newcommand{\de}{\delta}
\newcommand{\De}{\Delta}
\newcommand{\al}{\alpha}
\newcommand{\be}{\beta}
\newcommand{\ga}{\gamma}
\newcommand{\Ga}{\Gamma}
\newcommand{\la}{\lambda}
\newcommand{\La}{\Lambda}
\newcommand{\iy}{\infty}
\newcommand{\om}{\omega}
\newcommand{\Om}{\Omega}
\newcommand{\vfi}{\varphi}
\newcommand{\veps}{\varepsilon}
\newcommand{\kp}{\mbox{\rm\ae}}
\newcommand{\fv}{\mathfrak{M}}
\newcommand{\Res}{\mathop{\rm Res}\limits}
\newcommand{\br}[1]{{\left\langle {#1}\right\rangle}}
\newcommand{\bbr}[1]{{\Big\langle {#1}\Big\rangle}}
\newcommand{\te}{\theta}

\begin{document}

\thispagestyle{empty}

\begin{center}
{\large\bf
RECOVERING DIRAC SYSTEM WITH SINGULARITIES IN INTERIOR POINTS}\\[0.2cm]
{\bf O. Gorbunov and V.Yurko} \\[0.2cm]
\end{center}

{\bf Abstract.} We study the non-selfadjoint Dirac system on a finite interval
having non-integrable regular singularities in interior points with additional
matching conditions at these points. Properties of spectral characteristics are
established, and the inverse spectral problem is investigated. We provide a
constructive procedure for the solution of the inverse problem, and prove its
uniqueness. Moreover, necessary and sufficient conditions for the global
solvability of this nonlinear inverse problem are obtained.

Key words: differential systems, singularity, spectral analysis, inverse problems

AMS Classification: 34A55 34L40 34A36  47E05 \\

{\bf 1. Introduction. } Consider the boundary value problem
$L=L(Q_{\om}(x), Q(x), \al,\be)$ for the Dirac system on a finite
interval with $N$ regular singularities inside the interval:
$$
BY'+\Big(Q_\om(x)+Q(x)\Big)Y=\la Y,\quad 0<x<\pi,                             \eqno(1)
$$
$$
(\cos\alpha,\;\sin\alpha)Y(0)=(\cos\be,\;\sin\be)Y(\pi)=0,                    \eqno(2)
$$
where
$$
Y(x)=\left(\begin{array}{c} y_1(x) \\ y_2(x) \end{array}\right),\;
B=\left(\begin{array}{rc} 0&1\\-1&0\end{array}\right),\;
Q(x)=\left(\begin{array}{cr} q_1(x)&q_2(x)\\q_2(x)&-q_1(x)\end{array}\right),
$$
$Q_\om(x)=Q^\br{k}_\om(x)=\di\frac{\mu_k}{x-\ga_k}\left(\begin{array}{rr}
\sin2\eta_k&\cos2\eta_k\\ \cos2\eta_k&\sin2\eta_k\end{array}\right)$
for $x\in\om_{k+1/2}\cup\ga_{k+1/2},\;k=\overline{1,N}$.

Here $0<\ga_1<\ga_2<\ldots<\ga_N<\pi,\; \om_p=(\ga_{p},\,\ga_{p+1}),
\; \ga_{k+1/2}=(\ga_{k+1}+\ga_k)/2,\; k=\overline{1,\,N-1},\;
\ga_{1/2}=\ga_0=0,\; \ga_{N+1/2}=\ga_{N+1}=\pi$, $q_j(x)$ are
complex-valued functions, and $\mu_k$ are complex numbers. Let for
definiteness, $\alpha,\;\be,\;\eta_k\in[-\pi/2,\pi/2],$ $\mbox{Re}\,\mu_k>0,\;\mu_k+1/2\notin\mathbb{N}.$ Let $q_j(x)$ be
absolutely continuous on $[0,\pi]$ and $\di|q_j(x)|\prod\limits_{k=1}^N
|x-\ga_k|^{-2Re\mu_k}\in L(0,\pi)$. If $Q_{\om}(x), Q(x), \al,\be$
satisfy these conditions, we will say that  $L\in W.$

In this paper we establish properties of spectral characteristics
and investigate the inverse spectral problem of recovering $L$ from
the given spectral data. We provide a constructive procedure for the
solution of the inverse problem, and prove its uniqueness. Moreover,
necessary and sufficient conditions for the global solvability of this
nonlinear inverse problem are obtained.

Differential equations with singularities inside the interval play an
important role in various areas of mathematics as well as in applications.
Moreover, a wide class of differential equations with turning points can be
reduced to equations with singularities. For example, such problems appear
in electronics for constructing parameters of heterogeneous electronic lines
with desirable technical characteristics [1]-[3]. Boundary value problems
with discontinuities in an interior point appear in geophysical models for
oscillations of the Earth [4]. Differential equations with turning points
arise in various physical and technical problems; see [5] where further
references and links to applications can be found.
We also note that in different problems of natural sciences
we face different kind of matching conditions in singular points.

The case when a singular point lies at the endpoint of the interval
was investigated fairly completely for various classes of
differential equations in [6]-[10] and other works.
The presence of singularity inside the interval produces
essential qualitative modifications in the investigation (see [11]).

A few words on the structure of the paper. In section 2 properties
of spectral characteristics are studied. For this we use the results
from [12] where special fundamental systems of solutions are
constructed with prescribed analytic and asymptotic properties.
In section 3 we provide a constructive procedure for the
solution of the inverse problem, and prove its uniqueness.
Necessary and sufficient conditions for the global solvability of
the inverse problem are presented in section 4.

\medskip
{\bf 2. Properties of the spectrum. } System (1) has non-integrable
singularities at the points $\ga_k$, hence it is necessary to reguire
additional matching conditions for solutions on the intervals $\om_{k-1}$
and $\om_k$. We will do it as follows. It was shown in [12] that for
$x\in\om_{k-1}\cup\om_k$ there exist a fundamental system of solutions  $S^\br{k}(x,\la)=(S_1^\br{k}(x,\la),\; S_2^\br{k}(x,\la))$ such that
$$
S^\br{k}_1(x,\la)\sim(x-\ga_k)^{-\mu_k}
\left(\begin{array}{c}0\\c_{01}\end{array}\right),\;
S^\br{k}_2(x,\la)\sim(x-\ga_k)^{\mu_k}
\left(\begin{array}{c}c_{02}\\0\end{array}\right)
\;\mbox{for}\;x\rightarrow\ga_k.
$$
where $c_{01}c_{02}=1.$ Let  $Y(x,\la)=a_1(\la)S^\br{k}_1(x,\la)+
a_2(\la)S^\br{k}_2(x,\la)$ be a solution of system (1) for
$x\in\om_{k-1}$. Then we put by definition
$$
Y(x,\la)=a_1(\la)S^\br{k}_1(x,\la)A^\br{k}(\la)+
a_2(\la)S^\br{k}_2(x,\la)A^\br{k}(\la),
$$
for $x\in\om_k$, where $A^\br{k}(\la)$ is a fixed given transition matrix
for  $\ga_k$. For example, if $A^\br{k}(\la)=I$ ($I$ is the edentity matrix)
and $Q(x)$ is analytic at $\ga_k$, then this continuation of the solution
coincides with the analytic continuation through the upper half-plane
$Im x>0.$ If $A^\br{k}(\la)=
\left(\begin{array}{cc}e^{2i\pi\mu_k}&0\\0&e^{-2i\pi\mu_k}\end{array}\right)$,
then it corresponds to the analytic continuation through the lower half-plane
$Im x <0.$

Let $S(x,\la)=(S_1(x,\la),\;S_2(x,\la))$ be the fundamental matrix for
system (1) with the initial condition $S(0,\la)=I$  and with the above
mentioned matching conditions. For definiteness, everywhere below
$A^\br{k}(\la)=I,\,k=\overline{1,N}$. The construction of this fundamental
matrix can be described as follows. If $x\in\om_0\cup\om_1$, then we put
$S(x,\la)=S^\br{1}(x,\la)\Big(S^\br{1}(0,\la)\Big)^{-1}$;
moreover, if $x\in\om_1$, then $S(x,\la)=S^\br{2}(x,\la)C^\br{1}(\la).$
Fix $x_1\in\om_1$. Then  $S^\br{1}(x_1,\la)\Big(S^\br{1}(0,\la)\Big)^{-1}
=S^\br{2}(x_1,\la)C^\br{1}(\la)$, i.e.
$$
S(x,\la)=S^\br{2}(x,\la)\Big(S^\br{2}(x_1,\la)\Big)^{-1}
S^\br{1}(x_1,\la)\Big(S^\br{1}(0,\la)\Big)^{-1},\; x\in\om_1.
$$
Analogously, one gets for $x\in\om_k$:
$$
S(x,\la)=S^\br{k+1}(x,\la)
\left(\prod\limits_{j=1}^{k}\Big(S^\br{j+1}(x_j,\la)\Big)^{-1}
S^\br{j}(x_j,\la)\right) (S^\br{1}(0,\la))^{-1},\; x_j\in\om_j.                \eqno(3)
$$

{\bf Lemma 1. }{\it For $x\in\om_k$ and $|\la(x-\ga_k)|\ge1$,
$$
S(x,\la)=\frac{1}{2i}\Big(e^{i\la
x}\left[\begin{array}{rr}i&-1\\1&i\end{array}\right]_\br{k} + e^{-i\la x}
\left[\begin{array}{rr}i&1\\-1&i\end{array}\right]_\br{k}\Big) +
\sum\limits_{j=1}^k\sin\pi\mu_je^{-il\la(x-2\ga_j)+2il\eta_j}
\left[\begin{array}{rr}-i&l\\l&i\end{array}\right]_\br{k},
$$
$\left[\Big(a_{ij}\Big)_{i,j=1}^{n,m}\right]_\br{k} :=
\Bigg(a_{ij}+O\Big(|\la(x-\ga_k)|^{-\nu}\Big)\Bigg)_{i,j=1}^{n,m},\,\,
\nu=\min\{1,\;2{\rm Re}\mu_1,\;2{\rm Re}\mu_2,\;\ldots,2{\rm Re}\mu_N\},$\\
$l=\left\{\begin{array}{rl}1,&\arg\la\in\Pi_{-1}\cup\Pi_1,\\
-1,&\arg\la\in\Pi_0,\end{array}\right.,\,\,
\di\Pi_k=\left\{\la\;\Big|\;\arg\la\in\Big(\pi\frac{5k-3}{6-2k},
\pi\frac{5k+3}{6+2k}\Big]\right\},\,\,k=0,\;\pm1.$}

\smallskip
We prove the lemma by induction.
According to [12], the matrix $S^\br{k}(x,\la)$ can be represented by
$S^\br{k}(x,\la)=E^\br{k}(x,\la)\be^\br{k}(\la),$ where $E^\br{k}(x,\la)$
is the Birkhoff-type fundamental matrix, and $\be^\br{k}(\la)$ are Stockes
multipliers.

Let $x\in\om_0$. Then $S(x,\la)=S^\br{1}(x,\la)(S^\br{1}(0,\la))^{-1}$
and (see [12])
$$
S(x,\la)=\frac{1}{[2i]_\br{1}}
\left(\begin{array}{cc}e^{i\la x)}[i]_\br{1} + e^{-i\la x)}[i]_\br{1}&
e^{i\la x)}[-1]_\br{1} + e^{-i\la x)}[1]_\br{1}\\[2mm]
e^{i\la x)}[1]_\br{1} + e^{-i\la x)}[-1]_\br{1} &e^{i\la x)}[i]_\br{1}
+e^{-i\la x)}[i]_\br{1}\end{array}\right).
$$
Suppose that the assertion of the lemma is true for $x\in\om_{k-1}$.
Let us prove it for $x\in\om_k$. It follows from (3) that for  $x\in\om_k$,
$$
S(x,\la)=S^\br{k}(x,\la)(S^\br{k}(x_{k-1},\la))^{-1}S(x_{k-1},\la).           \eqno(4)
$$
We find the asymptotics for $S^\br{k}(x,\la)(S^\br{k}(x_{k-1},\la))^{-1},$
using the asymptotics from [12]. Denote $l^+=l^\br{k},\;m^+=m^\br{k}$
 for $x>\ga_k$, and $l^-=l^\br{k},\;m^-=m^\br{k}$ for $x<\ga_k$. One has
$$
S^\br{k}(x,\la)B\Big(S^\br{k}(x_{k-1},\la)\Big)^TB^T
$$
$$
=\left(e^{-i\la(x-\ga_k)+i\eta_k}\left[\begin{array}{rr}-i&-i\\1&1
\end{array}\right]_\br{k}+e^{i\la(x-\ga_k)-i\eta_k}\left[
\begin{array}{rr}i&-i\\1&-1\end{array}\right]_\br{k}H(e^{i\pi\mu_kl^+})\right)
$$
$$
\times H(e^{2i\pi\mu_km^+})H(\la^{-\mu_k})\be^\br{k}B\be^\br{k}H(\la^{-\mu_k})
H(e^{2i\pi\mu_km^-})
$$
$$
\times
\left(e^{-i\la(x_{k-1}-\ga_k)+i\eta_k}\left[\begin{array}{rr}-i&1\\-i&1
\end{array}\right]_\br{k}+e^{i\la(x_{k-1}-\ga_k)-i\eta_k}H(e^{i\pi\mu_kl^-})
\left[\begin{array}{rr}i&1\\-i&-1\end{array}\right]_\br{k}\right)B^T,
$$
where $H(z)=\left(\begin{array}{ll}z^{-1}&0\\0&z\end{array}\right)$,
and $T$ is the sign for the transposition.  Since
 $BH(z)B^T=H(z^{-1})$ and $\be^\br{k}B\be^\br{k}B^T
=\be^\br{k}_1\be^\br{k}_2$, it follows that
$$
S^\br{k}(x,\la)B\Big(S^\br{k}(x_{k-1},\la)\Big)^TB^T
=\be^\br{k}_1\be^\br{k}_2\Big(
e^{-i\la(x-\ga_k)+i\eta_k}\left[\begin{array}{rr}-i&-i\\1&1
\end{array}\right]_\br{k}H(e^{2i\pi\mu_k(m^+-m^-)})
$$
$$
+e^{2i\la(x-\ga_k)-i\eta_k}\left[\begin{array}{rr}i&-i\\1&-1
\end{array}\right]_\br{k}H(e^{i\pi\mu_k(l^++2m^+-2m^-)})\Big)
$$
$$
\times
\left(e^{-i\la(x_{k-1}-\ga_k)+i\eta_k}\left[\begin{array}{rr}1&i\\-1&-i
\end{array}\right]_\br{k}+e^{i\la(x_{k-1}-\ga_k)-i\eta_k}H(e^{-i\pi\mu_kl^-})
\left[\begin{array}{rr}-1&i\\-1&i\end{array}\right]_\br{k}\right).
$$
Taking the relation $\be^\br{k}_1\be^\br{k}_2=(4i\cos\pi\mu_k)^{-1}$
into account, we calculate
$$
S^\br{k}(x,\la)\Big(S^\br{k}(x_{k-1},\la)\Big)^{-1}
$$
$$
=\frac{1}{4i\cos\pi\mu_k}\; e^{-i\la(x+x_{k-1}-2\ga_k)+2i\eta_k}
\left[2i\sin\Big(2\pi\mu_k(m^--m^+)\Big)
\left(\begin{array}{rr}-i&1\\1&i\end{array}\right)\right]_\br{k}
$$
$$
+\frac{1}{4i\cos\pi\mu_k}\; e^{-i\la(x-x_{k-1})}
\left[2\cos\Big(\pi\mu_k(l^-+2m^--2m^+)\Big)
\left(\begin{array}{rr}i&1\\-1&i\end{array}\right)\right]_\br{k}
$$
$$
+\frac{1}{4i\cos\pi\mu_k}\; e^{i\la(x-x_{k-1})}
\left[2\cos\Big(\pi\mu_k(l^++2m^+-2m^-)\Big)
\left(\begin{array}{rr}i&-1\\1&i\end{array}\right)\right]_\br{k}
$$
$$
+\frac{1}{4i\cos\pi\mu_k}\; e^{i\la(x+x_{k-1}-2\ga_k)-2i\eta_k}
\left[2\sin\Big(\pi\mu_k(l^--l^++2m^--2m^+)\Big)
\left(\begin{array}{rr}-i&-1\\-1&i\end{array}\right) \right]_\br{k}.
$$
Consider three cases :\\
1) If $\la\in\Pi_1$, then $m^-=1,\,m^+=0,\,l^-=-1,\,l^+=1,$ and
$$
S^\br{k}(x,\la)\Big(S^\br{k}(x_{k-1},\la)\Big)^{-1}
$$
$$
= \sin(\pi\mu_k) e^{-i\la(x+x_{k-1}-2\ga_k)+2i\eta_k}
\left[\begin{array}{rr}-i&1\\1&i\end{array}\right]_\br{k}+ \frac{1}{2i}\;
e^{-i\la(x-x_{k-1})} \left[\begin{array}{rr}i&1\\-1&i\end{array}
\right]_\br{k}
$$
$$
+\frac{1}{2i}\; e^{i\la(x-x_{k-1})} \left[
\begin{array}{rr}i&-1\\1&i\end{array}\right]_\br{k}+
e^{i\la(x+x_{k-1}-2\ga_k)-2i\eta_k}
\left[\begin{array}{rr}0&0\\0&0\end{array}\right]_\br{k}.
$$
2) If $\la\in\Pi_{-1}$, then $m^-=0,\,m^+=-1,\,l^-=-1,\,l^+=1,$ and
$$
S^\br{k}(x,\la)\Big(S^\br{k}(x_{k-1},\la)\Big)^{-1}
$$
$$
= \sin(\pi\mu_k) e^{-i\la(x+x_{k-1}-2\ga_k)+2i\eta_k}
\left[\begin{array}{rr}-i&1\\1&i\end{array}\right]_\br{k}+ \frac{1}{2i}\;
e^{-i\la(x-x_{k-1})} \left[\begin{array}{rr}i&1\\-1&i\end{array}\right]_\br{k}
$$
$$
+\frac{1}{2i}\; e^{i\la(x-x_{k-1})} \left[
\begin{array}{rr}i&-1\\1&i\end{array}\right]_\br{k}+
e^{i\la(x+x_{k-1}-2\ga_k)-2i\eta_k}
\left[\begin{array}{rr}0&0\\0&0\end{array}\right]_\br{k}.
$$
3) If $\la\in\Pi_0$, then $m^-=0,\,m^+=0,\,l^-=1,\,l^+=-1,$ and
$$
S^\br{k}(x,\la)\Big(S^\br{k}(x_{k-1},\la)\Big)^{-1}
= e^{-i\la(x+x_{k-1}-2\ga_k)+2i\eta_k}
\left[\begin{array}{rr}0&0\\0&0\end{array}\right]_\br{k}+ \frac{1}{2i}\;
e^{-i\la(x-x_{k-1})} \left[\begin{array}{rr}i&1\\-1&i\end{array}
\right]_\br{k}
$$
$$
+\frac{1}{2i}\; e^{i\la(x-x_{k-1})} \left[
\begin{array}{rr}i&-1\\1&i\end{array}\right]_\br{k}+
\sin(\pi\mu_k) e^{i\la(x+x_{k-1}-2\ga_k)-2i\eta_k}
\left[\begin{array}{rr}-i&-1\\-1&i\end{array}\right]_\br{k}.
$$
Since $x_{k-1}<\ga_k,\,x>\ga_k$, it follows that  $x+x_{k-1}-2\ga_k=
x-x_{k-1}+2(x_{k-1}-\ga_k)<x-x_{k-1}$, $x+x_{k-1}-2\ga_k=
x_{k-1}-x+2(x-\ga_k)>x_{k-1}-x$, and the exponentials
$e^{\pm i\la(x+x_{k-1}-2\ga_k)}$ grow not faster than
$e^{\pm i\la(x-x_{k-1})}$. Thus,
$$
S^\br{k}(x,\la)\Big(S^\br{k}(x_{k-1},\la)\Big)^{-1}=
 \frac{1}{2i}\; e^{i\la(x-x_{k-1})}
\left[\begin{array}{rr}i&-1\\1&i\end{array}\right]_\br{k}
$$
$$
+\frac{1}{2i}\; e^{-i\la(x-x_{k-1})} \left[
\begin{array}{rr}i&1\\-1&i\end{array}\right]_\br{k}+
\sin(\pi\mu_k) e^{-li\la(x+x_{k-1}-2\ga_k)+2li\eta_k}
\left[\begin{array}{rr}-i&l\\l&i\end{array}\right]_\br{k}.
$$
Substituting this asymptotics into (4), we get
$$
S(x,\la)=\left( \frac{1}{2i}\; e^{i\la(x-x_{k-1})}
\left[\begin{array}{rr}i&-1\\1&i\end{array} \right]_\br{k}+ \frac{1}{2i}\;
e^{-i\la(x-x_{k-1})} \left[
\begin{array}{rr}i&1\\-1&i\end{array}\right]_\br{k}\right.
$$
$$
\left. +\sin(\pi\mu_k) e^{-li\la(x+x_{k-1}-2\ga_k)+2li\eta_k}
\left[\begin{array}{rr}-i&l\\l&i\end{array}\right]_\br{k} \right)
\left(\frac{1}{2i}\;e^{i\la
x_{k-1}}\left[\begin{array}{rr}i&-1\\1&i\end{array}\right]_\br{k-1}\right.
$$
$$
\left. +\frac{1}{2i}\;e^{-i\la x_{k-1}}
\left[\begin{array}{rr}i&1\\-1&i\end{array}\right]_\br{k-1} +
\sum\limits_{j=1}^{k-1}\sin\pi\mu_je^{-il\la(x_{k-1}-2\ga_j)+2il\eta_j}
\left[\begin{array}{rr}-i&l\\l&i\end{array}\right]_\br{k-1}\right).
$$
Since $0<x_{k-1}<x,$ it follows that $0<2x_{k-1}<2x,\;
-x<2x_{k-1}-x<x,$ and $e^{\pm i\la(2x_{k-1}-x)}$ grow not faster than
$e^{\pm i\la x}$. Therefore
$$
S(x,\la)=\frac{1}{2i}\;e^{i\la
x}\left[\begin{array}{rr}i&-1\\1&i\end{array}\right]_\br{k}+
\frac{1}{2i}\;e^{-i\la
x}\left[\begin{array}{rr}i&1\\-1&i\end{array}\right]_\br{k}
$$
$$
+\frac{1}{2i}\sum\limits_{j=1}^{k-1}\sin(\pi\mu_j)
e^{i\la(x-x_{k-1})-li\la(x_{k-1}-2\ga_j)+2li\eta_j}
\left[(1-l)\left(\begin{array}{rr}1&-i\\-i&-1\end{array}\right)\right]_\br{k}
$$
$$
+\frac{1}{2i}\sum\limits_{j=1}^{k-1}\sin(\pi\mu_j)
e^{-i\la(x-x_{k-1})-li\la(x_{k-1}-2\ga_j)+2li\eta_j}
\left[(1+l)\left(\begin{array}{rr}1&i\\i&-1\end{array}\right)\right]_\br{k}
$$
$$
+\frac{1}{2i}\sin(\pi\mu_k) e^{i\la x_{k-1}-li\la(x+x_{k-1}-2\ga_k)+2li\eta_k}
\left[(1+l)\left(\begin{array}{rr}1&i\\i&-1\end{array}\right)\right]_\br{k}
$$
$$
+\frac{1}{2i}\sin(\pi\mu_k) e^{-i\la x_{k-1}-li\la(x+x_{k-1}-2\ga_k)+2li\eta_k}
\left[(1-l)\left(\begin{array}{rr}1&-i\\-i&-1\end{array}\right)\right]_\br{k}
$$
$$
+\sum\limits_{j=1}^{k-1}\sin(\pi\mu_k)\sin(\pi\mu_j)
e^{-li\la(x+2x_{k-1}-2\ga_k-2\ga_j)+4li\eta_j}
\left[\begin{array}{rr}0&0\\0&0\end{array}\right]_\br{k}.
$$
Let $l=-1.$ Then
$$
S(x,\la)=\frac{1}{2i}\;e^{i\la x}
\left[\begin{array}{rr}i&-1\\1&i\end{array}\right]_\br{k}+
\frac{1}{2i}\;e^{-i\la x}
\left[\begin{array}{rr}i&1\\-1&i\end{array}\right]_\br{k}
$$
$$
+\sum\limits_{j=1}^{k-1}\sin(\pi\mu_j) e^{i\la x-2i\la\ga_j-2i\eta_j}
\left[\begin{array}{rr}-i&-1\\-1&i\end{array}\right]_\br{k}
+\sum\limits_{j=1}^{k-1}\sin(\pi\mu_j) e^{-i\la x+2i\la
x_{k-1}-2i\la\ga_j-2i\eta_j}
\left[\begin{array}{rr}0&0\\0&0\end{array}\right]_\br{k}
$$
$$
+ e^{2i\la x_{k-1}+i\la x-2i\la\ga_k-2i\eta_k}
\left[\begin{array}{rr}0&0\\0&0\end{array}\right]_\br{k}+ \sin(\pi\mu_k)
e^{i\la x-2i\la\ga_k-2i\eta_k}
\left[\begin{array}{rr}1&-i\\-i&-1\end{array}\right]_\br{k}
$$
$$
+\sum\limits_{j=1}^{k-1}\sin(\pi\mu_k)\sin(\pi\mu_j)
e^{i\la(x+2x_{k-1}-2\ga_k-2\ga_j)-4i\eta_j}
\left[\begin{array}{rr}0&0\\0&0\end{array}\right]_\br{k}.
$$
This yields
$$
S(x,\la)=\frac{1}{2i}\Big(e^{i\la
x}\left[\begin{array}{rr}i&-1\\1&i\end{array}\right]_\br{k}
+e^{-i\la
x}\left[\begin{array}{rr}i&1\\-1&i\end{array}\right]_\br{k}\Big)
+\sum\limits_{j=1}^k\sin(\pi\mu_j) e^{i\la x-2i\la\ga_j-2i\eta_j}
\left[\begin{array}{rr}-i&-1\\-1&i\end{array}\right]_\br{k}.
$$
The case $l=1$ is treated similarly. Lemma 1 is proved.

\smallskip
The following assertion is proved analogously.

\smallskip
{\bf Lemma 2. }{\it For $x\in\om_k$ and $|\la(x-\ga_k)|\ge1$
$$
\frac{\partial}{\partial\la} S(x,\la)=\frac{x}{2i}
\left(e^{i\la x}\left[\begin{array}{rr}-1&-i\\i&-1\end{array}\right]_\br{k}
+e^{-i\la x}\left[\begin{array}{rr}1&-i\\i&1\end{array}\right]_\br{k}\right)
$$
$$
+\sum\limits_{j=1}^k(x-2\ga_j)\sin\pi\mu_je^{-il\la(x-2\ga_j)+2il\eta_j}
\left[\begin{array}{rr}-l&-i\\-i&l\end{array}\right]_\br{k}.
$$}

{\bf Definition. } A function $Y(x,\la)$ is called the solution of system (1),
if there exist constants $C_1(\la),\,C_2(\la)$ such that
$Y(x,\la)=C_1(\la)S_1(x,\la)+C_2(\la)S_2(x,\la),$
$x\in(0,\pi)\setminus\bigcup\limits_{k=1}^N\{\ga_k\}.$

We introduce the functions
$$
\vfi(x,\la)=\Big(\vfi_1(x,\la),\;\vfi_2(x,\la)\Big)=S(x,\la)V(\al),\;
V(\al)=\Big(V_1(\al),V_2(\al)\Big)=\left(\begin{array}{rr}\cos\al&-\sin\al\\
\sin\al&\cos\al\end{array}\right),
$$
$$
\De(\la)=\left(\begin{array}{cc}\De_{11}(\la)&\De_{12}(\la)\\
\De_{21}(\la)&\De_{22}(\la)\end{array}\right)=V^T(\be)S(\pi,\la)V(\al),
$$
$$
\psi(x,\la)=\Big(\psi_1(x,\la),\;\psi_2(x,\la)\Big)=S(x,\la)S^{-1}(\pi,\la)V(\be).
$$
Clearly,  $\vfi(x,\la),\,\psi(x,\la)$ are fundamental matrices for
system (1). Denote $\br{Y,Z}:=Y^TBZ.$ If $Y(x,\la), Z(x,\la)$ are solutions
of system (1), then  $\br{Y(x,\la),Z(x,\la)}:=\det\{Y(x,\la),Z(x,\la)\}$ is
their Wroskian. Obviously,
$$
\br{\psi_2(x,\la),\vfi_2(x,\la)}=-\De_{12}(\la).                             \eqno(5)
$$
A number $\la_0$ is called an eigenvalue of problem (1)-(2), if
there exist constants $A_1,\,A_2$ ($|A_1|+|A_2|>0$) such that
the function $A_1S_1(x,\la_0)+A_2S_2(x,\la_0)$ satisfies
the boundary conditions (2).

\smallskip
{\bf Lemma 3. }{\it Zeros of $\De_{12}(\la)$ coincide with the eigenvalues
of the boundary value problem (1)-(2). If $\la_0$ is an eigenvalue, then
$\vfi(x,\la_0)$ and $\psi(x,\la_0)$ are eigenfunctions, and
$\psi(x,\la_0) = b_0\vfi(x,\la_0)$.}

{\it Proof.} 1) Let $\la_0$ be a zero of $\De_{12}(\la),$ i.e.
$V_1^T(\be)S(\pi,\la_0)V_2(\al)=0.$ Therefore,
$\vfi_2(x,\la_0)=S(x,\la_0)V_2(\al)$ is an eigenfunction,
and $\la_0$ is an eigenvalue. It follows from (5) that
$\vfi_2(x,\la_0)$ and $\psi_2(x,\la_0)$ are linear dependent.

2) Let $\la_0$ be an eigenvalue, and let $Y_0(x)$ be the corresponding
eigenfunction. Since $\vfi_1(x,\la),\,\vfi_2(x,\la)$ form
a fundamrntal system of solutions, it follows that \\
$Y_0(x)=D_1\vfi_1(x,\la_0)+D_2\vfi_2(x,\la_0).$
Substituting this relation into the first boundary condition, we obtain
$D_1V_1^T(\al)V_1(\al)+D_2V_1^T(\al)V_2(\al)=0,$ hence
$D_1=0$. Using the second boundary condition, we find
$D_2V_1^T(\be)\vfi_2(\pi,\la_0)=0.$
Since $Y_0(x)\not\equiv0,$ one has $D_2\neq0,$ i.e.\\
$V_1^T(\be)\vfi_2(\pi,\la_0)=0.$ Lemma 3 is proved.

\smallskip
We note that the functions $\De_{jk}(\la),\;j,k=1,2,$ are the characteristic
functions for the boundary value problems $L_{jk}$ foe system (1) with
boundary conditions. $V^T_{3-k}(\al)Y(0)=V^T_j(\be)Y(\pi)=0.$ Denote
$$
\Phi(x,\la)=\Big(\Phi_1(x,\la),\,\Phi_2(x,\la)\Big),\;
\Phi_1(x,\la)=-\frac{1}{\De_{12}(\la)}\psi_2(x,\la),\,
\Phi_2(x,\la)=\vfi_2(x,\la).
$$
It follows from (5) that $\det\Phi(x,\la)\equiv1.$ The functions
$\Phi_1(x,\la),\,\Phi_2(x,\la)$ are called the Weyl solutions, and
the matrix $\fv(\la):=V^T_2(\la)\Phi_1(0,\la)$ is called the Weyl
matrix for the problem (1)-(2).

\smallskip
{\bf Lemma 4. } {\it The following relations hold
$$
\Phi(x,\la)=\vfi(x,\la)M(\la),\;\mbox{where }\;
M(\la)=\left(\begin{array}{cc}1&0\\ \fv(\la)&1\end{array}\right),\;
\fv(\la)=-\frac{\De_{11}(\la)}{\De_{12}(\la)}.
$$}

Only formula for $\Phi_1(x,\la)$ is needed to be proved. Let \\
$\Phi_1(x,\la)=D_1(\la)\vfi_1(x,\la)+D_2(\la)\vfi_2(x,\la).$ Then
$$
-(\De_{12}(\la))^{-1}\psi_2(0,\la)=D_1(\la)V_1(\al)+D_2(\la)V_2(\al).      \eqno(6)
$$
Multiplying (6) by $V_1^T(\al),$ we infer
$-(\De_{12}(\la))^{-1}V_1^T(\al)\psi_2(0,\la)=D_1(\la).$ Since
$$
V_1^T(\al)\psi_2(0,\la)=V_1^T(\al)B^TS^T(\pi,\la)BV_2(\be),
$$
it follows that  $V_1^T(\al)\psi_2(0,\la)=-V_1^T(\be)S(\pi,\la)V_2(\al)
=-\De_{12}(\la),$ i.e. $D_1(\la)=1$. Multiplying (6) by $V_2^T(\al),$
we find $D_2(\la)=V_2^T(\al)\Phi_1(0,\la)=\fv(\la).$
Taking the relation $V_2^T(\al)\psi_2(0,\la)=\De_{11}(\la)$
into account, we obtain the assertion of the lemma.

Thus, $\fv(\la)$ is a meromorphic function; its poles coincide with the
eigenvalues of $L,$ and its zeros coincide with the eigenvalues of $L_{11}$.

\smallskip
{\bf Lemma 5. } {\it For $x\in\om_k$ and
$|\la(x-\ga_k)|\ge1,\,|\la(x-\ga_{k+1})|\ge1,$ one has
$$
\vfi(x,\la)=\frac{1}{2i}\left(e^{i\la
x+i\al}\left[\begin{array}{rr}i&-1\\1&i\end{array}\right]_\br{k}+e^{-i\la
x-i\al} \left[\begin{array}{rr}i&1\\-1&i\end{array}\right]_\br{k}\right)
$$
$$
+\sum_{j=1}^k \sin\pi\mu_je^{-il\la(x-2\ga_j)+2il\eta_j-il\al}
\left[\begin{array}{rr}-i&l\\l&i\end{array}\right]_\br{k},                \eqno(7)
$$
$$
\psi(x,\la)=-\frac{1}{2i}\left(e^{i\la(\pi-x)-i\be}
\left[\begin{array}{rr}i&1\\-1&i\end{array}\right]_\br{k+1} +
e^{-i\la(\pi-x)+i\be}
\left[\begin{array}{rr}i&-1\\1&i\end{array}\right]_\br{k+1}\right)
$$
$$
-\sum_{j=k+1}^N \sin\pi\mu_je^{-il\la(x+\pi-2\ga_j)+2il\eta_j-il\be}
\left[\begin{array}{rr}-i&-l\\-l&i\end{array}\right]_\br{k+1}.            \eqno(8)
$$}

Indeed, since $\vfi(x,\la)=S(x,\la)V(\al)$, relation (7) follows from Lemma 1.
To prove (8) we make the substitution $x\to \pi-x$ and repeat the arguments.

Taking the relation $\De(\la)=V(\be)\vfi(\pi,\la)$ into account we arrive at
the following assertion.

\smallskip
{\bf Corollary 1. } {\it For the characteristic function
$\De_{12}(\la),$ the following asymptotics holds
$$
\De_{12}(\la)=\frac{1}{2i}e^{-i(\la\pi+\al-\be)}[1]-
\frac{1}{2i}e^{i(\la\pi+\al-\be)}[1]+\sum_{j=1}^N
\sin\pi\mu_j e^{-il\la(\pi-2\ga_j)+il(2\eta_j-\al-\be)}[l],             \eqno(9)
$$
where $\left[\left(a_{kj}\right)_{k,j=1}^{n,m}\right]:=
\Big(a_{kj}+O\left(|\la|^{-\nu}\right)\Big)_{k,j=1}^{n,m}$
for $|\la|\to\iy$.}

\smallskip
By the well-known method (see for example [14-15]) one
obtains the following properties:

1. $\De_{12}(\la)=O(e^{\pi|Im\la|})$.

2. All eigenvalues $\la_k,\;k\in\mathbb{Z}$ of the problem
(1)-(2) lie in the strip $|Im\la|\le h.$

3. Let $N_a$ be a number of eigenvalues in the rectange\\
$\Big\{\la\,|\,Re\la\in[a,a+1),\;|Im\la|\le h\Big\}.$ Then
$N_a$ is uniformly bounded.

4. Denote $G_\de=\{\la\,:\,|\la-\la_k|\ge\de\;\forall\, k\}$. Then
$|\De_{12}(\la)|\ge C_\de e^{\pi|Im\la|}$ for $\la\in G_{\de}.$

5. For sufficiently small $\de,$ there exists a sequence
$R_n\to\iy$ such that the circles $\Ga_n=\Big\{\la\,:\,|\la|
=R_n\Big\}$ lie in $G_\de.$

6. Let $\left\{\la^0_k\right\}_{k=-\iy}^\iy$ be zeros of the function
$$
\De^0_{12}(\la)=\frac{1}{2i}e^{-i(\la\pi+\al-\be)}-
\frac{1}{2i}e^{i(\la\pi+\al-\be)} + l\sum_{j=1}^N
\sin\pi\mu_j e^{-il\la(\pi-2\ga_j)+il(2\eta_j-\al-\be)}.                \eqno(10)
$$
Then $\la_k=\la^0_k+O(|\la^0_k|^{-\nu}).$

\smallskip
For simplicity, we confine ourselves to the case when all
eigenvalues of $L$ are simple, i.e. the function $\Delta_{12}(\la)$
has only simple zeros. In particular, it is always true for the
self-adjoint case. Denote $a_k:=\Res_{\la=\la_k}\,\fv(\la).$
The data $\{a_k,\;\la_k\}_{k=-\iy}^{+\iy}$ are called the
spectral data for $L.$ The inverse problem is formulated
as follows.

\smallskip
{\bf Inverse problem 1.} Given $\{a_k,\;\la_k\}_{k=-\iy}^{+\iy},$
construct $L,$ i.e. $Q(x),\;Q_\om(x),\;\al,\;\be.$

\smallskip
In Sections 3-4 we give an algorithim for the global solution
of this nonlinear inverse problem and provide necessary and
sufficient conditions for its solvability.

\smallskip
{\bf Lemma 6. }{\it Let $\fv^0(\la)$ be the Weyl function for
the problem $L^0$ of the form (1)-(2) but with the zero potential
$Q(x)\equiv 0.$ Then
$$
\fv(\la)=\fv^0(\la)+\sum_{k=-\iy}^{+\iy}\left(\frac{a_k}{\la-\la_k}-
\frac{a^0_k}{\la-\la^0_k}\right),\quad \di\sum_{k=-\iy}^{+\iy}=
\di\lim_{n\to\iy}\di\sum_{|\la_k|<R_n,\,|\la^0_k|<R_n}.
$$

Proof. } Consider the integral $J_n(\la)=\di\frac{1}{2\pi i}\int_{\Ga_n}
\frac{\fv(\xi)-\fv^0(\xi)}{\xi-\la}d\xi$, $\la\in\mbox{int}\,\Ga_n$.
Using lemmas 4-5 and Corollary 1, we obtain
$\fv(\xi)-\fv^0(\xi)=O(|\xi|^{-\nu})$ for $\xi\in G_\de,$
and consequently, $J_n(\la)\to0$ as $n\to\iy.$
On the othe hand, by residue's theorem,
$$
J_n(\la)=\Res_{\xi=\la}\frac{\fv(\xi)-\fv^0(\xi)}{\xi-\la}+
\sum_{|\la|<R_n,\,|\la^0_k|<R_n}\left(\Res_{\xi=\la_k}
\frac{\fv(\xi)}{\xi-\la}-
\Res_{\xi=\la^0_k}\frac{\fv^0(\xi)}{\xi-\la}\right),
$$
hence
$$
J_n(\la)=\fv(\la)-\fv^0(\la)+
\sum_{|\la|<R_n,\,|\la^0_k|<R_n}\left(\frac{a_k}{\la_k-\la}-
\frac{a^0_k}{\la^0_k-\la}\right).
$$
If $n\to\iy,$ we arrive at the assertion of the lemma.

\smallskip
Together with $L$ we consider a boundary value problem $\widetilde L$
of the same form (1)-(2) but with different $\widetilde Q(x),\;\widetilde Q_\om(x)$,
$\widetilde\al,\;\widetilde\be$. We agree that if a certain symbol $v$
denotes an object related to $L,$ then $\widetilde v$ will denote an
analogous object related to $\widetilde L.$

\smallskip
{\bf Lemma 7. }{\it If $\la_k=\widetilde\la_k$ for all $k,$ then
$\De_{12}(\la)\equiv\widetilde\De_{12}(\la).$

Proof. } The functions $\De_{12}(\la)$ and $\widetilde\De_{12}(\la)$
are entire in $\la$ of exponential type. Using Hadamard's factorization theorem,
we get $\De_{12}(\la)=e^{a\la+b}\widetilde\De_{12}(\la).$
Let us show that $a=0,\;b=0.$ In view of (9),
$$
\frac{1}{2i}e^{-i(\la\pi+\al-\be)}[1]-\frac{1}{2i}e^{i(\la\pi+\al-\be)}[1]
+\sum_{j=1}^N \sin\pi\mu_j e^{-il\la(\pi-2\ga_j)+il(2\eta_j-\al-\be)}[l]=
$$
$$
\frac{1}{2i}e^{-i(\la\pi+\widetilde\al-\widetilde\be)+a\la+b}[1]-
\frac{1}{2i}e^{i(\la\pi+\widetilde\al-\widetilde\be)+a\la+b}[1]+
\sum_{j=1}^N \sin\pi\widetilde\mu_je^{-il\la(\pi-2\widetilde\ga_j)
+il(2\widetilde\eta_j-\widetilde\al-\widetilde\be)+a\la+b}[l].                   \eqno(11)
$$
Let $\la=\sigma+i\tau.$ If $\tau=0$ and $\sigma\to+\iy,$ then the
right-hand side in (11) is bounded; hence $Re\,a\le0$; for $\tau=0$ and
$\sigma\to-\iy,$ we get $Re\,a\ge0,$ i.e. $Re\,a=0.$ Furthermore,
the right-hand side in (11) is $O(e^{-Im\la\tau}),$ but the left-hand side
is $O(e^{-Im\la\tau -Im\,aIm\la})$ for $\tau\le0.$ For $\tau\to-\iy$
we have $Im\,a\le0.$ If $\tau\ge0,$ then it follows from (11) that
$O(e^{\pi\tau})=O(e^{\pi\tau-Im\,a\tau}).$ This means that
$Im\,a\ge0,$ i.e. $Im\,a=0.$ Thus, $a=0.$ Similarly, one gets that
 $b=0.$ Lemma is proved.

\smallskip
{\bf Corollary 2. }{\it If $\la_k=\widetilde\la_k$ for all $k,$
then $\De^0_{12}(\la)\equiv\widetilde\De^0_{12}(\la),$ i.e.
$\al-\be=\widetilde\al-\widetilde\be,\;\ga_k=\widetilde\ga_k,$
$\sin\pi\mu_ke^{il(2\eta_k-\al-\be)}=
\sin\pi\widetilde\mu_ke^{il(2\widetilde\eta_k-\widetilde\al-
\widetilde\be)}$. Here $\De^0_{12}(\la)$ is defined by (10).}

\smallskip
{\bf Lemma 8. }{\it If
$\al-\widetilde\al=\be-\widetilde\be=\widetilde\eta_k-\eta_k,\;
\mu_k=\widetilde\mu_k,\;\ga_k=\widetilde\ga_k,\;k=\overline{1,N}$
and $Q(x)=\widetilde Q(x)V^2(\widetilde\al-\al),$ then
$\fv(\la)=\widetilde\fv(\la)$.

Proof. } Denote $\de:=\al-\widetilde\al=\be-\widetilde\be=\widetilde\eta_k-\eta_k$.
Let us show that if $Y(x,\la)$ is a solution of (1), then $\widetilde Y(x,\la)=
V(-\de)Y(x,\la)$ is a solution of $\widetilde{(1)}$. Indeed, substituting
$V(\de)\widetilde Y(x,\la)$ into (1), we obtain
$$
BV(\de)\widetilde Y'(x,\la)+\Big(Q(x)+Q_\om(x)\Big)V(\de)\widetilde
Y(x,\la)=\la V(\de)\widetilde Y(x,\la).
$$
Multiplying by $V^T(\de)=V(-\de)$ and taking the relation
$V^T(\de)Q(x)=Q(x)V(\de)$ into account, we get
$$
B\widetilde Y'(x,\la)+\Big(Q(x)+Q_\om(x)\Big)
V^2(\de)\widetilde Y(x,\la)=\la\widetilde Y(x,\la).
$$
One has $V(-\de)S(0,\la)V(\de)=I.$ Since the Cauchy problem has the unique
solution, we infer $\widetilde S(x,\la)=V(-\de)S(x,\la)V(\de).$ Then
$\widetilde\De(\la)=V^T(\widetilde\be)V(-\de)S(\pi,\la)V(\de)
V(\widetilde\al)$ or
$$
\widetilde\De(\la)=V^T(\widetilde\be+\de)S(\pi,\la)V(\de+\widetilde\al).
$$
This yields $\widetilde\De_{jk}(\la)=\De_{jk}(\la).$ Lemma is proved.

\medskip
{\bf 3. Solution of the inverse problem. } Let us first prove
the uniqueness theorem.

\smallskip
{\bf Theorem 1. }{\it If $\fv(\la)=\widetilde\fv(\la),$ then
$\al-\widetilde\al=\be-\widetilde\be=\widetilde\eta_k-\eta_k,\;
\mu_k=\widetilde\mu_k, \;\ga_k=\widetilde\ga_k,\;k=\overline{1,N}$
and $Q(x)=\widetilde Q(x)V^2(\widetilde\al-\al).$

Proof. } By virtue of lemma 8, it is sufficient to prove the theorem
for the case $\widetilde\be=\be=0.$ Consider the function $P(x,\la)=
\Phi(x,\la)\widetilde\Phi^{-1}(x,\la).$

Since $\fv(\la)=\widetilde\fv(\la),$ it follows thaty these functions
have the same poles. In view of lemma 7, one gets $\De_{12}(\la)=
\widetilde\De_{12}(\la).$ By corollary 2, $\widetilde\al=\al,\;
\widetilde\ga_k=\ga_k,\; \sin\pi\mu_ke^{il(2\eta_k-\al)}=
\sin\pi\widetilde\mu_ke^{il(2\widetilde\eta_k-\widetilde\al)}.$ This yields
$$
\vfi(x,\la)-\widetilde\vfi(x,\la)=O(e^{|Im\la|x}|\la|^{-\nu}),\;\;\;
\psi(x,\la)-\widetilde\psi(x,\la)=O(e^{|Im\la|(\pi-x)}|\la|^{-\nu}).     \eqno(12)
$$
Since $\Phi_1(x,\la)=-(\De_{12}(\la))^{-1}\psi_2(x,\la),$ it follows that
$$
\Phi_1(x,\la)=O(e^{-x|Im\la|}),\;\;\la\in G_\de.                         \eqno(13)
$$
Taking (12) into account, we infer
$$
\Phi_1(x,\la)-\widetilde\Phi_1(x,\la)=
O(e^{-x|Im\la|}|\la|^{-\nu}),\;\;\la\in G_\de.                           \eqno(14)
$$
Obviously, $P(x,\la)-I=
(\Phi(x,\la)-\widetilde\Phi(x,\la))B\widetilde\Phi(x,\la)B^T.$
Using (12)-(14), we obtain for  $\la\in G_\de:$
$$
P(x,\la)-I=|\la|^{-\nu}
\left(\begin{array}{rr}O(e^{-|Im\la|x})&O(e^{|Im\la|x})\\
O(e^{-|Im\la|x})&O(e^{|Im\la|x})\end{array}\right)
\left(\begin{array}{rr}O(e^{|Im\la|x})&O(e^{|Im\la|x})\\
O(e^{-|Im\la|x})&O(e^{-|Im\la|x})\end{array}\right)
=O(|\la|^{-\nu}).                                                       \eqno(15)
$$
Since $\Phi(x,\la)=\vfi(x,\la)M(\la),$ we get $P(x,\la)=
\vfi(x,\la)M(\la)\widetilde M^{-1}(\la)\widetilde\vfi^{-1}(x,\la),$
or $P(x,\la)=\vfi(x,\la)\widetilde\vfi^{-1}(x,\la).$ Therefore,
$P(x,\la)$ is entire in $\la.$ Using (15), maximum modulus principle
and Liouville's theorem, we conclude that $P(x,\la)=I,$ i.e.
$\widetilde\Phi(x,\la)=\Phi(x,\la).$ Then
$Q(x)+Q_\om(x)=\widetilde Q(x)+\widetilde Q_{\widetilde\om}(x),$
and consequently,
$Q(x)=\widetilde Q(x),$ $Q_\om(x)=\widetilde Q_{\widetilde\om}(x).$
Theorem is proved.

\smallskip
{\bf Corollary 3. }{\it If $a_k=\widetilde a_k,\;\la_k=
\widetilde\la_k$ for all $k,$ then $L=\widetilde L.$}

\smallskip
{\bf Corollary 4. }
{\it If $\la^\br{11}_k=\widetilde\la^\br{11}_k,\;\la_k=\widetilde\la_k$
for all $k,$ then $L=\widetilde L.$ Here
$\{\la^\br{11}_k\}$ and $\{\widetilde\la^\br{11}_k\}$ are eigenvalues
of $L_{11}$ and $\widetilde L_{11}$, respectively.}

\smallskip
Indeed, according to lemma 7, $\De_{12}(\la)=\widetilde\De_{12}(\la).$
Analogously, we obtain $\De_{11}(\la)=\widetilde\De_{11}(\la).$
By lemma 4, $\fv(\la)=\widetilde\fv(\la).$

\smallskip
Let us now go on to constructing the solution of the nonlinear inverse
problem 1. The central role here is played by the so-called main
equation of the inverse problem, which is a linear equation in the
corresponding Banach space. Let us derive the main equation.

Let the problem $L$ with a simple spectrum be given. We choose
a model boundary value problem $\widetilde L$ with a simple
spectrum such that $\om=\tilde\om,$
$Q_\om(x)=\widetilde Q_{\om}(x)$ and
$$
\La := \sum_{k=-\iy}^{+\iy}|\widetilde a_k|\xi_k <\iy,\;\;
\xi_k:=|\la_k-\widetilde\la_k|+ |\widetilde a^{-1}_ka_k-1|.            \eqno(16)
$$
For definiteness, we assume that
$\al=\widetilde\al=0.$ Then $\be=\widetilde\be.$
Denote
$\Om_\veps:=\{x\,:\,x\in(0,\pi),\,|x-\ga_k|\ge\veps,\;
k=\overline{1,N}\}$, $\la_{k0}=\la_k,\,\la_{k1}
=\widetilde\la_k,\;a_{k0}=a_k,\,a_{k1}=\widetilde a_k$,
$$
\widetilde D^\br{l}(x,\la,\te):=\frac{\br{\widetilde\Phi_l(x,\la),
\; \widetilde\vfi_2(x,\te)}}{\la-\te},\;l=1,2,\quad
\widetilde D^\br{2}_{kj}(x,\la)=\widetilde D^\br{2}(x,\la,\la_{kj}),
$$
$$
\widetilde P_{ni,kj}(x)=\widetilde D^\br{2}(x,\la_{ni},\la_{kj})a_{kj},
\; \vfi_{2,kj}(x)=\vfi_2(x,\la_{kj}),\;
\widetilde\vfi_{2,kj}(x)=\widetilde\vfi_2(x,\la_{kj}),
$$
where $\br{Y,Z}:=\det(Y,Z)=Y^TBZ.$ Analogously we define
$D^\br{l}(x,\la,\te),$ $D^\br{2}_{kj}(x,\la)$ and $P_{ni,kj}(x).$

\smallskip
{\bf Lemma 9. }{\it For $x\in\Om_\veps$ ans $\la$ on compact sets,
$$
|\widetilde\vfi^{(m)}_{2,kj}(x)|\le C(1+|\la^0_k|)^m,\;\;
|\widetilde\vfi^{(m)}_{2,k1}(x)-\widetilde\vfi^{(m)}_{2,k0}(x)|
\le C\xi_k(1+|\la^0_k|)^m,\;m=0,1,                                     \eqno(17)
$$
$$
\left.\begin{array}{c} |\widetilde D^\br{2}_{kj}(x,\la)|
\le \di\frac{C}{1+|\la-\la^0_k|},\;\;
|\widetilde D^\br{2}_{k0}(x,\la)a_{k0}-
\widetilde D^\br{2}_{k1}(x,\la)a_{k1}|
\le\di\frac{C|a_{k1}|\xi_k}{1+|\la-\la^0_k|},\\[5mm]
|(\widetilde D^\br{2}_{kj}(x,\la))'|\le C,\;\;
|(\widetilde D^\br{2}_{k0}(x,\la)a_{k0}-
\widetilde D^\br{2}_{k1}(x,\la)a_{k1})'|\le C|a_{k1}|\xi_k.
\end{array}\right\}.                                                    \eqno(18)
$$
The same estimates are valid for $\vfi_{2,kj}(x),\;D^\br{2}_{kj}(x,\la)$.}

\smallskip
In order to prove the lemma, we need the following generalization of
Schwarz's lemma:\\
{\it Let the function $f(z)$ be analytic inside the circle $|z-z_0|\le R$
and continuous in the whole circle. Moreover, $|f(z)|\le C$ on the boundary,
and $f(z_0)=0.$ Then $|f(z)|\le C|z-z_0|/R$ in the circle $|z-z_0|\le R$.}

1) It follows from (7) that
$$
|\widetilde\vfi_2(x,\la)|\le C e^{|Im\la|x},\;\; x\in\Om_\veps.        \eqno(19)
$$
The eigenvalues lie in the strip $|Im\la|\le\max\{h,\,\widetilde h\}$;
it follows from (19) that
$|\widetilde\vfi^{(m)}_{2,kj}(x)|\le C(1+|\la_{kj}|)^m.$
Using (10), we obtain the first estimate in (17) for $m=0.$

Applying Schwarz's lemma, we find
$|\widetilde\vfi_2(x,\la)-\widetilde\vfi_{2,k1}(x)|\le
Ce^{|Im\la|x}|\la-\la_{k1}|.$ Hence the second estimate
in (17) holds (17) for $m=0.$ For $m=1,$
the arguments are similar.

2) Since $\widetilde D^\br{2}(x,\la,\te)=(\la-\te)^{-1}
(\widetilde\vfi_2(x,\la))^T B\widetilde\vfi_2(x,\te),$
it follows from (19) for $\la\neq\te$, $|\la|\le R,\;|\te|\le R$
that $|\widetilde D(x,\la,\te)|\le C|\la-\te|^{-1}.$
If $\la=\te,$ then $\widetilde D(x,\la,\la)=(\widetilde\vfi_2(x,\la))^T
B\dot{\widetilde\vfi}_2(x,\la),$ where $\dot{\widetilde\vfi}_2(x,\la)=
\di\frac{\partial}{\partial\la}\widetilde\vfi_2(x,\la).$ Using lemma 2,
we obtain $|\dot{\widetilde\vfi}_2(x,\la)|\le Cxe^{|Im\la|x}$
for $x\in\Om_\veps.$
Then $|\widetilde D(x,\la,\la)|\le C$ for $|\la|\le R.$
Thus,
$$
|\widetilde D(x,\la,\te)|\le\frac{C}{1+|\la-\te|},
\;\;x\in\Om_\veps,\;\;|\la|\le R,\;\;|\te|\le R.                    \eqno(20)
$$
Furthermore,
$\br{\widetilde\Phi_j(x,\la),\widetilde\vfi_2(x,\te)}'=
(\widetilde\Phi^T_j(x,\la))'B\widetilde\vfi_2(x,\te)+
\widetilde\Phi^T_j(x,\la)B\widetilde\vfi\,'_2(x,\te).$ Then
$$
\br{\widetilde\Phi_j(x,\la),\widetilde\vfi_2(x,\te)}'=
-(B\widetilde\Phi'_j(x,\la))^T\widetilde\vfi_2(x,\te)+
\widetilde\Phi^T_j(x,\la)B\widetilde\vfi\,'_2(x,\te).
$$
Since $\widetilde\Phi_j, \widetilde\vfi_2$ are solutions
of the system, it follows that\\
$\br{\widetilde\Phi_j(x,\la),\widetilde\vfi_2(x,\te)}'= (\te-\la)
\widetilde\Phi^T_j(x,\la)\widetilde\vfi_2(x,\te).$
This yields
$$
(\widetilde D^\br{j}(x,\la,\te))'
=-\widetilde\Phi^T_j(x,\la)\widetilde\vfi_2(x,\te).                \eqno(21)
$$
Taking (21) and (19) into account, we arrive at
the third estimate in (18).

Using Schwarz's lemma and (20), we infer
$|\widetilde D^\br{2}_{k0}(x,\la)-\widetilde D^\br{2}_{k1}(x,\la)|
\le\di\frac{C|\la_{k0}-\la_{k1}|}{1+|\la-\la^0_k|}.$
Since
$$
|\widetilde D^\br{2}_{k0}(x,\la)a_{k0}-\widetilde D^\br{2}_{k1}(x,\la)a_{k1}|
\le |\widetilde D^\br{2}_{k0}(x,\la)(a_{k0}-a_{k1})|+
|(\widetilde D^\br{2}_{k0}(x,\la)-\widetilde D^\br{2}_{k1}(x,\la))a_{k1}|,
$$
one gets the second estimate in (18). Other estimates are obtained
analogously. Lemma is proved.

\smallskip
Similarly one can prove the following assertion.

\smallskip
{\bf Lemma 10. }{\it For $x\in\Om_\veps$ and $\la$
on compact sets,
$$
|\widetilde P_{ni,kj}(x)|\le\frac{C|a_{k1}|}{1+|\la^0_n-\la^0_k|},
\;\; |\widetilde P'_{ni,kj}(x)|\le C|a_{k1}|,
$$
$$
|\widetilde P_{ni,k1}(x)-\widetilde P_{ni,k0}(x)|\le
\frac{C|a_{k1}|\xi_k}{1+|\la^0_n-\la^0_k|},\;\;
|\widetilde P'_{ni,k1}(x)-\widetilde P'_{ni,k0}(x)|\le C|a_{k1}|\xi_k,
$$
$$
|\widetilde P_{n1,kj}(x)-\widetilde P_{n0,kj}(x)|\le
\frac{C|a_{k1}|\xi_n}{1+|\la^0_n-\la^0_k|},\;\;
|\widetilde P'_{n1,kj}(x)-\widetilde P'_{n0,kj}(x)|\le C|a_{k1}|\xi_n,
$$
$$
|\widetilde P_{n1,k1}(x)-\widetilde P_{n1,k0}(x)-\widetilde P_{n0,k1}(x)
+\widetilde P_{n0,k0}(x)| \le\frac{C|a_{k1}|\xi_k\xi_n}{1+|\la^0_n-\la^0_k|},
$$
$$
|\widetilde P'_{n1,k1}(x)-\widetilde P'_{n1,k0}(x)-\widetilde P'_{n0,k1}(x)
+\widetilde P'_{n0,k0}(x)| \le C|a_{k1}|\xi_k\xi_n.
$$
Moreover, if $\la\in G_\de=\Big\{\la\,:
\,|\la-\widetilde\la_k|\ge\de,\;k\in\mathbb{Z}\Big\},$ then
$$
|\widetilde D^\br{1}_{kj}(x,\la)|\le\frac{C_\de}{|\la-\la_{kj}|},\;\;
|(\widetilde D^\br{1}_{kj}(x,\la))'|\le C_\de,
$$
$$
|\widetilde D^\br{1}_{k0}(x,\la)a_{k0}-\widetilde D^\br{1}_{k1}(x,\la)
a_{k1}|\le C_\de|a_{k1}|\xi_k\left(\frac{1}{|\la-\la_{k0}|}+
\frac{1}{|\la-\la_{k1}|}\right),
$$
$$
|(\widetilde D^\br{1}_{k0}(x,\la)a_{k0}-\widetilde D^\br{1}_{k1}(x,\la)
a_{k1})'|\le C_\de|a_{k1}|\xi_k,
$$
where $C$ and $C_\de$ depend on $\veps.$ The same estimates
are valid for $D^\br{l}_{kj}(x,\la),\;P_{ni,kj}(x)$.}

\smallskip
{\bf Lemma 11. }{\it The following relations hold
$$
\Phi_j(x,\la)=\widetilde\Phi_j(x,\la)+\sum_{k=-\iy}^{+\iy}\Big(
\widetilde D^\br{j}_{k0}(x,\la)a_{k0}\vfi_{2,k0}(x)-
\widetilde D^\br{j}_{k1}(x,\la)a_{k1}\vfi_{2,k1}(x)\Big),\; j=1,2,    \eqno(22)
$$
the series converge absolutely and uniformly for $x\in\Om_\veps$ and
$\la$ on compact sets without the spectra of $L$ and $\widetilde L.$

\smallskip
Proof. }
Consider the function $P(x,\la)=\Phi(x,\la)\widetilde\Phi^{-1}(x,\la).$
Denote
$$
J_n(x,\la)=\di\frac{1}{2\pi i}di\int_{\Ga_n}\di\frac{1}{\xi-\la}
\Big(P(x,\xi)-I\Big)d\xi,\quad \Ga_n:=\Big\{\la\,:\,|\la|=R_n\Big\}.
$$
The functions $\Phi(x,\la)$ and $\widetilde\Phi(x,\la)$ have the same
main term in the asymptotics. Therefore, for a fixed $x\neq\ga_k$,
one has $P(x,\xi)-I=O(|\xi|^{-\nu}),$ and $J_n(x,\la)\to0$ as $n\to\iy$
uniformly in $\la$ on the compact sets. Integration on  $\Ga_n$
is divided into integration on the contours
$\Ga^\br{1}_n=\Ga^\br{3}_n\bigcup\Ga^\br{5}_n$,
$\Ga^\br{2}_n=\Ga^\br{4}_n\bigcup\Ga^\br{5}_n$
(with counterclockwise circuit), where
  $\Ga^\br{3}_n=\{\la\,:\,|Im\la|\le h\}\bigcap\Ga_n$,
$\Ga^\br{4}_n=\Ga_n\setminus\Ga^\br{3}_n=
\{\la\,:\,|Im\la|>h\}\bigcap\Ga_n$,
$\Ga^\br{5}_n=\{\la\,:\,|Im\la|=h\}\bigcap\mbox{int}\,\Ga_n$.
Let $\la\in\mbox{int}\,\Ga^\br{2}_n$.
By the Cauchy integral formula,
$$
\frac{1}{2\pi i}\int_{\Ga^\br{2}_n}\frac{1}{\xi-\la}
\Big(P(x,\xi)-I\Big)d\xi = P(x,\la)-I.
$$
Clearly,
$$
\frac{1}{2\pi i}\int_{\Ga^\br{1}_n}\frac{1}{\xi-\la}Id\xi = 0.
$$
Then
$$
P(x,\la)=I+\frac{1}{2\pi i}
\int_{\Ga^\br{1}_n}\frac{1}{\la-\xi}P(x,\xi)d\xi-J_n(x,\la).           \eqno(23)
$$
Since $\Phi(x,\la)=\vfi(x,\la)M(\la),$ it follows that
$$
P(x,\xi)=\vfi(x,\xi)M(\xi)\widetilde M^{-1}(\xi)\widetilde\vfi^{-1}(x,\xi)
$$
$$
=\vfi(x,\xi)\widetilde\vfi^{-1}(x,\xi)-
(\fv(\xi)-\widetilde\fv(\xi))\vfi(x,\xi)B_\br{1}\widetilde\vfi^{-1}(x,\xi),
\; B_\br{1}=\left(\begin{array}{rr}0&0\\-1&0\end{array}\right).
$$
The function $\vfi(x,\xi)\widetilde\vfi^{-1}(x,\xi)$ is entire in $\xi.$
Therefore,
$$
\frac{1}{2\pi i}\int_{\Ga^\br{1}_n}\vfi(x,\xi)\widetilde\vfi^{-1}(x,\xi)
\frac{d\xi}{\la-\xi}=0,
$$
since $\la$ is outside $\Ga^\br{1}_n$. Thus, it follows from (23) that
$$
P(x,\la)=I-\frac{1}{2\pi i}\int_{\Ga^\br{1}_n}
(\fv(\xi)-\widetilde\fv(\xi))\vfi(x,\xi)B_\br{1}\widetilde\vfi^{-1}(x,\xi)
\frac{d\xi}{\la-\xi}-J_n(x,\la).
$$
One has $\Phi(x,\la)=P(x,\la)\widetilde\Phi(x,\la),$ hence
$\Phi_j(x,\la)=P(x,\la)\widetilde\Phi_j(x,\la).$ Then
$$
\Phi_j(x,\la)=\widetilde\Phi_j(x,\la)-\frac{1}{2\pi i}\int_{\Ga^\br{1}_n}
(\fv(\xi)-\widetilde\fv(\xi))\vfi(x,\xi)B_\br{1}\widetilde\vfi^{-1}(x,\xi)
\widetilde\Phi_j(x,\la)\frac{d\xi}{\la-\xi}+\veps_n(x,\la),
$$
and $\veps_n(x,\la)\to0$ as $n\to\iy$ uniformly for $x\in\Om_\veps$.
Furthermore,
$$
B_\br{1}\widetilde\vfi^{-1}(x,\xi) \widetilde\Phi_j(x,\la)=
\Big(\widetilde\vfi_{12}(x,\xi)\widetilde\Phi_{2j}(x,\la)-
\widetilde\vfi_{22}(x,\xi)\widetilde\Phi_{1j}(x,\la)\Big)
\left(\begin{array}{c}0\\1\end{array}\right),
$$
and consequently,
$$
\Phi_j(x,\la)=\widetilde\Phi_j(x,\la)-\frac{1}{2\pi i}\int_{\Ga^\br{1}_n}
(\fv(\xi)-\widetilde\fv(\xi))\br{\widetilde\Phi_j(x,\la),
\widetilde\vfi_2(x,\xi)}\vfi_2(x,\xi)\frac{d\xi}{\la-\xi}+\veps_n(x,\la).
$$
Calculating the integral by residue's theorem and taking $n\to\iy,$
we arrive at (22). Lemma is proved.

\smallskip
Consider (22) for $j=2$ and $\la=\la_{ni}$:
$$
\widetilde\vfi_{m2,ni}(x)=\vfi_{m2,ni}(x)-
\sum_{k=-\iy}^{+\iy}\Big(\widetilde P_{ni,k0}(x)\vfi_{m2,k0}(x)-
\widetilde P_{ni,k1}(x)\vfi_{m2,k1}(x)\Big),\,\,m=1,2,                   \eqno(24)
$$
where $\vfi_{2,kj}(x)=
\left(\begin{array}{c}\vfi_{12,kj}(x)\\\vfi_{22,kj}(x)\end{array}\right)$.
The last relation is not convenient for our purpose, since the series
converges  only "with brackets".  We trasform (24) as follows:
$$
\widetilde\vfi_{m2,n0}(x)-\widetilde\vfi_{m2,n1}(x)=
\vfi_{m2,n0}(x)-\vfi_{m2,n1}(x)
-\sum_{k=-\iy}^{+\iy}\Big((\widetilde P_{n0,k0}(x)-
\widetilde P_{n1,k0}(x))(\vfi_{m2,k0}(x)-\vfi_{m2,k1}(x))
$$
$$
+(\widetilde P_{n0,k0}(x)-\widetilde P_{n1,k0}(x)-
\widetilde P_{n0,k1}(x)+\widetilde P_{n1,k1}(x))\vfi_{m2,k1}(x)\Big),
$$
$$
\widetilde\vfi_{m2,n1}(x)=\vfi_{m2,n1}(x)-\sum_{k=-\iy}^{+\iy}
\Big(\widetilde P_{n1,k0}(x)(\vfi_{m2,k0}(x)-\vfi_{m2,k1}(x))+
(\widetilde P_{n1,k0}(x)-\widetilde P_{n1,k1}(x))\vfi_{m2,k1}(x)\Big).
$$
Denote
$$
\Psi^\br{m}_{n0}(x)=\chi_n\Big(\vfi_{m2,n0}(x)-\vfi_{m2,n1}(x)\Big),
\; \chi_n=\left\{\begin{array}{rl}0,&\xi_n=0,\\\xi^{-1}_n,&\xi_n\neq0,
\end{array}\right.\; \Psi^\br{m}_{n1}(x)=\vfi_{m2,n1}(x),
$$
$$
\widetilde H_{n0,k0}(x)=(\widetilde P_{n0,k0}(x)-
\widetilde P_{n1,k0}(x))\chi_n\xi_k,
$$
$$
\widetilde H_{n0,k1}(x)=(\widetilde P_{n0,k0}(x)-\widetilde P_{n1,k0}(x)
-\widetilde P_{n0,k1}(x)+\widetilde P_{n1,k1}(x))\chi_n,
$$
$$
\widetilde H_{n1,k0}(x)=\widetilde P_{n1,k0}(x)\xi_k,\;
\widetilde H_{n1,k1}(x)=\widetilde P_{n1,k0}(x)-\widetilde P_{n1,k1}(x).
$$
Then
$$
\widetilde\Psi^\br{m}_{ni}(x)=\Psi^\br{m}_{ni}(x)-\sum_{k=-\iy}^{+\iy}
\Big(\widetilde H_{ni,k0}(x)\Psi^\br{m}_{k0}(x)+
\widetilde H_{ni,k1}(x)\Psi^\br{m}_{k1}(x)\Big).                          \eqno(25)
$$
Using lemmas 9 and 10, we obtain the estimates
$$
\left.\begin{array}{c}
|\widetilde\Psi^\br{m}_{ni}(x)|\le C,\;
|(\widetilde\Psi^\br{m}_{ni}(x))'|\le C(1+|\la^0_k|),\\[2mm]
|\widetilde H_{ni,kj}(x)|
\le\di\frac{C|a_{k1}|\xi_k}{1+|\la^0_n-\la^0_k|},\,;
|\widetilde H'_{ni,kj}(x)|\le C|a_{k1}|\xi_k.
\end{array}\right\}                                                       \eqno(26)
$$
The same estimates are valid for $\Psi^\br{m}_{ni}(x),
\,\,H_{ni,kj}(x)$. Denote
$$
\Psi^\br{m}(x)=\left(\begin{array}{c}\Psi^\br{m}_{n0}(x)\\
\Psi^\br{m}_{n1}(x)\end{array}\right)_{n=-\iy}^{+\iy}=
\Big(\ldots,\Psi^\br{m}_{-1,1}(x),\,\;
\Psi^\br{m}_{00}(x),\!\Psi^\br{m}_{01}(x),\,\;
\Psi^\br{m}_{10}(x),\ldots\Big)^T.
$$
Similarly we define the block-matrix
$$
\widetilde H(x)=\left(\begin{array}{cc}\widetilde H_{n0,k0}(x)&
\widetilde H_{n0,k1}(x)\\ \widetilde H_{n1,k0}(x)&
\widetilde H_{n1,k1}(x)\end{array}\right)_{n,k=-\iy}^{+\iy}.
$$
Then we rewrite (25) as follows
$$
\widetilde\Psi^\br{m}(x)=(I-\widetilde H(x))\Psi^\br{m}(x),\;m=1,2,      \eqno(27)
$$
where $I$ is the identity operator. It follows from (26) that
$\Psi^\br{m}(x),\;\widetilde\Psi^\br{m}(x)\in {\bf m}$ for
each fixed $x\neq\ga_k,\;k=\overline{1,N},$ where
${\bf m}$ is the Banach space of bounded sequences.
The operator $\widetilde H(x),$ acting from ${\bf m}$ to ${\bf m},$
is a linear bounded operator, and
$$
\|\widetilde H(x)\|_{m\to m}\le C\sup_n\sum_{k=-\iy}^{+\iy}
\frac{|a_{k1}|\xi_k}{1+|\la_n^0-\la^0_k|}\le
C\sum_{k=-\iy}^{+\iy}|a_{k1}|\xi_k<\iy.
$$
For each fixed $x,$ relation (27) is a linear equation in ${\bf m}$
with respect to $\Psi^\br{m}(x).$ This equation is called
the main equation of the inverse problem.

\smallskip
{\bf Lemma 12. }{\it The following relation holds
$$
Q(x)=\widetilde Q(x)+B\kp(x)-\kp(x)B,                                  \eqno(28)
$$
where
$$
\kp(x)=\sum_{k=-\iy}^{+\iy}\Big(a_{k0}\widetilde\vfi_{2,k0}(x)
\vfi^T_{2,k0}(x)-a_{k1}\widetilde\vfi_{2,k1}(x)\vfi^T_{2,k1}(x)\Big),  \eqno(29)
$$
and the series converges uniformly for $x\in\Om_\veps.$

\smallskip
Proof. } Differentiating (22), we calculate
$$
\Phi'_j(x,\la)=\widetilde\Phi'_j(x,\la)+\sum_{k=-\iy}^{+\iy}\Big
((\widetilde D^\br{j}_{k0}(x,\la))'a_{k0}\vfi_{2,k0}(x)+
\widetilde D^\br{j}_{k0}(x,\la)a_{k0}\vfi\,'_{2,k0}(x)
$$
$$
-(\widetilde D^\br{j}_{k1}(x,\la))'a_{k1}\vfi_{2,k1}(x)-
\widetilde D^\br{j}_{k1}(x,\la)a_{k1}\vfi\,'_{2,k1}(x)\Big).
$$
Multiplying this relation by $B$ and using (21), we obtain
$$
\Big(\la I-Q(x)-Q_\om(x)\Big)\Phi_j(x,\la)=
\Big(\la I-\widetilde Q(x)-\widetilde Q_\om(x)\Big)\widetilde\Phi_j(x,\la)
$$
$$
+\sum_{k=-\iy}^{+\iy}
\Big(-\widetilde\Phi^T_j(x,\la)\widetilde\vfi_{2,k0}(x)a_{k0}B\vfi_{2,k0}(x)
+\widetilde D^\br{j}_{k0}(x,\la)a_{k0}\Big(\la_{k0}I-Q(x)-Q_\om(x)\Big)
\vfi_{2,k0}(x)
$$
$$
+\widetilde\Phi^T_j(x,\la)\widetilde\vfi_{2,k1}(x)a_{k1}B\vfi_{2,k1}(x)-
\widetilde D^\br{j}_{k1}(x,\la)a_{k1}\Big(\la_{k1}I-Q(x)-Q_\om(x)\Big)
\vfi_{2,k1}(x)\Big),
$$
and consequently,
$$
(Q(x)-\widetilde Q(x))\widetilde\Phi_j(x,\la)+\sum_{k=-\iy}^{+\iy}
(-\widetilde\Phi^T_j(x,\la)\widetilde\vfi_{2,k0}(x)a_{k0}B\vfi_{2,k0}(x)
+\widetilde D^\br{j}_{k0}(x,\la)a_{k0}(\la_{k0}-\la)\vfi_{2,k0}(x)
$$
$$
+\widetilde\Phi^T_j(x,\la)\widetilde\vfi_{2,k1}(x)a_{k1}B\vfi_{2,k1}(x)
-\widetilde D^\br{j}_{k1}(x,\la)a_{k1}(\la_{k1}-\la)\vfi_{2,k1}(x))=0.
$$
Since $\widetilde D^\br{j}_{ni}(x,\la)=
\di\frac{\widetilde\Phi^T_j(x,\la)B\widetilde\vfi_{2,ni}(x)}{\la-\la_{ni}}$,
it follows that
$$
(Q(x)-\widetilde Q(x))\widetilde\Phi_j(x,\la)+\sum_{k=-\iy}^{+\iy}
\Big(-\{\widetilde\Phi^T_j(x,\la)\widetilde\vfi_{2,k0}(x)a_{k0}\}B\vfi_{2,k0}(x)
$$
$$
-\{\widetilde\Phi^T_j(x,\la)B\widetilde\vfi_{2,k0}(x)a_{k0}\}\vfi_{2,k0}(x)
+\{\widetilde\Phi^T_j(x,\la)\widetilde\vfi_{2,k1}(x)a_{k1}\}B\vfi_{2,k1}(x)
$$
$$
+\{\widetilde\Phi^T_j(x,\la)B\widetilde\vfi_{2,k1}(x)a_{k1}\}\vfi_{2,k1}(x)\Big)=0.
$$
The matrices $Q(x)$ and $\widetilde Q(x)$ are symmetrical. Then
$$
\widetilde\Phi^T(x,\la)\Big\{Q(x)-\widetilde Q(x)+\sum_{k=-\iy}^{+\iy}
\Big(\Big(a_{k0}\widetilde\vfi_{2,k0}(x)\vfi^T_{2,k0}(x)-
a_{k1}\widetilde\vfi_{2,k1}(x)\vfi^T_{2,k1}(x)\Big)B
$$
$$
-B\Big(a_{k0}\widetilde\vfi_{2,k0}(x)\vfi^T_{2,k0}(x)-
a_{k1}\widetilde\vfi_{2,k1}(x)\vfi^T_{2,k1}(x)\Big)\Big)\Big\}=0.
$$
Multiplying by $(\widetilde\Phi^T(x,\la))^{-1},$ we arrive at (28).
It follows from the estimate
$$
|a_{k0}\widetilde\vfi_{2,k0}(x)\vfi^T_{2,k0}(x)-
a_{k1}\widetilde\vfi_{2,k1}(x)\vfi^T_{2,k1}(x)|\le
|\widetilde\vfi_{2,k0}(x)\vfi^T_{2,k0}(x)-
\widetilde\vfi_{2,k1}(x)\vfi^T_{2,k1}(x)|\cdot|a_{k1}|
$$
$$
+|\widetilde\vfi_{2,k1}(x)\vfi^T_{2,k1}(x)|
\cdot|a_{k0}-a_{k1}|\le C|a_{k1}|\xi_k
$$
that the series in (29) converges uniformly. Lemma is proved.

Let us now study the solvability of the main equation. For this
purpose we need the following assertion.

{\bf Lemma 13. }{\it The following relation holds
$$
D^\br{2}(x,\la,\te)=\widetilde D^\br{2}(x,\la,\te)+\sum_{k=-\iy}^{+\iy}
(\widetilde D^\br{2}_{k0}(x,\la)D^\br{2}_{k0}(x,\te)a_{k0}
-\widetilde D^\br{2}_{k1}(x,\la) D^\br{2}_{k1}(x,\te)a_{k1}),          \eqno(30)
$$
and the series converges uniformly for $x\in\Om_\veps$ and
 $\la$ on compact sets.

\smallskip
Proof. } According to (23) we have for $\la,\,\te\in\Ga^\br{2}_n$:
$$
P(x,\la)-P(x,\te)=\frac{1}{2\pi i}\int_{\Ga^\br{1}_n}
\left(\frac{1}{\la-\xi}-\frac{1}{\te-\xi}\right)P(x,\xi)d\xi+J_n(x,\la,\te),
$$
where $J_n(x,\la,\te)\to0$ as $n\to\iy$ uniformly for $x\in\Om_\veps$
and $\la,\;\te$ on compact sets. Therefore,
$$
\frac{1}{\la-\te}\Big(P^T(x,\la)-P^T(x,\te)\Big)=
\frac{1}{2\pi i}\int_{\Ga^\br{1}_n}
\frac{1}{(\la-\xi)(\xi-\te)}P^T(x,\xi)d\xi+J^1_n(x,\la,\te).           \eqno(31)
$$
Since $P(x,\xi)=\Phi(x,\xi)\widetilde\Phi^{-1}(x,\xi)=
-\Phi(x,\la)B\widetilde\Phi^T(x,\xi)B,$ it follows that
$$
\widetilde\vfi^T_2(x,\la)P^T(x,\xi)B\vfi_2(x,\te)=
-\widetilde\vfi^T_2(x,\la)B\widetilde\Phi(x,\xi)
B\Phi^T(x,\xi)B\vfi_2(x,\te).
$$
One has $\br{y,z}=y^TBz,$ and consequently,
$$
\widetilde\vfi^T_2(x,\la)P^T(x,\xi)B\vfi_2(x,\te)=
\br{\widetilde\vfi_2(x,\la),\widetilde\vfi_2(x,\xi)}
\bbr{\Phi_1(x,\xi),\vfi_2(x,\te)}
$$
$$
-\br{\widetilde\vfi_2(x,\la),\widetilde\Phi_1(x,\xi)}
\br{\vfi_2(x,\xi),\vfi_2(x,\te)}.                                     \eqno(32)
$$
Since $\bbr{\Phi_1(x,\la),\vfi_2(x,\la)}\equiv1$,
$\bbr{\vfi_2(x,\la),\vfi_2(x,\la)}\equiv0,$ we infer
$$
\widetilde\vfi^T_2(x,\la)P^T(x,\la)B\vfi_2(x,\te) =
\bbr{\vfi_2(x,\la),\vfi_2(x,\te)},
$$
$$
\widetilde\vfi^T_2(x,\la)P^T(x,\te)B\vfi_2(x,\te) =
\bbr{\widetilde\vfi_2(x,\la),\widetilde\vfi_2(x,\te)}.
$$
Multiplying (31) by $\widetilde\vfi^T_2(x,\la)$ from the left,
and by $B\vfi_2(x,\te)$ from the right, and using (32),
we calculate
$$
\frac{\bbr{\vfi_2(x,\la),\vfi_2(x,\te)}}{\la-\te}-
\frac{\bbr{\widetilde\vfi_2(x,\la),\widetilde\vfi_2(x,\te)}}{\la-\te}
=\frac{1}{2\pi i}\int_{\Ga^\br{1}_n}\Big(
\frac{\br{\widetilde\vfi_2(x,\la),\widetilde\vfi_2(x,\xi)}
\bbr{\Phi_1(x,\xi),\vfi_2(x,\te)}}{(\la-\xi)(\xi-\te)}
$$
$$
-\frac{\br{\widetilde\vfi_2(x,\la),\widetilde\Phi_1(x,\xi)}
\br{\vfi_2(x,\xi),\vfi_2(x,\te)}}{(\la-\xi)(\xi-\te)}\Big)d\xi+J^2_n(x,\la,\te).
$$
By lemma 4, $\Phi_1(x,\xi)=\vfi_1(x,\xi)+\fv(\xi)\vfi_2(x,\xi).$
This yields
$$
\frac{\bbr{\vfi_2(x,\la),\vfi_2(x,\te)}}{\la-\te}-
\frac{\bbr{\widetilde\vfi_2(x,\la),\widetilde\vfi_2(x,\te)}}{\la-\te}
$$
$$
=\frac{1}{2\pi i}\int_{\Ga^\br{1}_n}
\frac{\bbr{\widetilde\vfi_2(x,\la),\widetilde\vfi_2(x,\xi)}
\bbr{\vfi_2(x,\xi),\vfi_2(x,\te)}}{(\la-\xi)(\xi-\te)}
\Big(\fv(\xi)-\widetilde\fv(\xi)\Big)d\xi + J^2_n(x,\la,\te),
$$
since the integrals from analytic functions are equal to zero.
Calculating the integral by residue's theorem and taking $n\to\iy,$
we arrive at  (30) firstly for $|\la|\ge h,$ and by analytic
continuation for all $\la.$ Lemma is proved.

\smallskip
Taking $\la=\la_{ni},\;\;\te=\la_{lj}$ in (30) and multiplying
by $a_{lj}$, we obtain
$$
P_{ni,lj}(x)-\widetilde P_{ni,lj}(x)-\sum_{k=-\iy}^{+\iy}
\Big(\widetilde P_{ni,k0}(x)P_{k0,lj}(x)-
\widetilde P_{ni,k1}(x)P_{k1,lj}(x)\Big)=0.                         \eqno(33)
$$
Symmetrically, one has
$$
P_{lj,ni}(x)-\widetilde P_{lj,ni}(x)-\sum_{k=-\iy}^{+\iy}
\Big(P_{lj,k0}(x)\widetilde P_{k0,ni}(x)-
P_{lj,k1}(x)\widetilde P_{k1,ni}(x)\Big)=0.                         \eqno(34)
$$
It follows from (33)-(34) that
$$
H_{ni,lj}(x)-\widetilde H_{ni,lj}(x)-\sum_{k=-\iy}^{+\iy}
\Big(\widetilde H_{ni,k0}(x)H_{k0,lj}(x)-
\widetilde H_{ni,k1}(x)H_{k1,lj}(x)\Big)=0,                        \eqno(35)
$$
$$
H_{ni,lj}(x)-\widetilde H_{ni,lj}(x)-\sum_{k=-\iy}^{+\iy}
\Big(H_{ni,k0}(x)\widetilde H_{k0,lj}(x)-
H_{ni,k1}(x)\widetilde H_{k1,lj}(x)\Big)=0.                        \eqno(36)
$$
We rewrite relations (35) and (36) in the matrix form
$$
H(x)-\widetilde H(x)-\widetilde H(x)H(x)=0,\;\;\;\;
H(x)-\widetilde H(x)-H(x)\widetilde H(x)=0
$$
or
$(I-\widetilde H(x))(I+H(x))=I,\; (I+H(x))(I-\widetilde H(x))=I.$
Thus, we have proved the following assertion.

\smallskip
{\bf Theorem 2. }{\it For each fixed $x$ ($x\ne \ga_k,\;
k=\overline{1,N}$), the linear bounded operator
$I-\widetilde H(x),$ acting from ${\bf m}$ to ${\bf m},$
has the unique inverse operator, and the main equation
(27) is uniquely solvable in ${\bf m}$.}

\smallskip
The solution of Inverse problem 1 can be constructed
by the following algorithm.

\smallskip
{\bf Algorithm 1. } Given the spectral data
$\{\la_k,\;a_k\}_{k=-\iy}^{+\iy}$ of the problem $L.$

1) Choose a model boundary value problem $\tilde L,$ for example,
with the zero potential.

2) Construct $\widetilde\Psi^\br{m}(x)$ and $\widetilde H(x).$

3) Solving the linear main equation (27), find $\Psi^\br{m}(x),$ and
then calculate $\vfi_{2,kj}(x).$

4) Construct $Q(x)$ by (28), and $\al=\widetilde\al,\;\;\be=\widetilde\be.$

\medskip
{\bf 4. Necessary and sufficient conditions for the solvability of
the inverse problem. }

\smallskip
{\bf Theorem 3. }{\it For numbers $\{\la_k,\;a_k\}_{k=-\iy}^{+\iy},$
$a_k\neq0$, $\la_k\neq\la_n$, $(k\neq n),$ to be the spectral data for
a certain problem $L\in W,$ it is necessary and sufficient that the
following conditions hold

1) (Asymptotics): There exists $\widetilde L\in W$ such that (16) holds;

2) (Condition S): For each fixed $x\neq\ga_k,\;k=\overline{1,N},$
the linear bounded operator $I-\tilde H(x)$ has the unique inverse operator;

3) $(B\kp(x)-\kp(x)B)|x-\ga_k|^{-2Re\mu_k}\in L(w_{k+1/2}),$ where
$\kp(x)$ is constructed by (29).

Under these conditions the potential $Q(x)$ is constructed by (28) and
$\al=\widetilde\al$, $\be=\widetilde\be$.}

\smallskip
The necessity part of the theorem was proved above. Let us prove the
sufficiency. Let numbers $\{\la_k,\;a_k\}_{k=-\iy}^{+\iy}$ be given
such that $a_k\neq0$ and $\la_k\neq\la_n$, $(k\neq n).$ Let
$\widetilde L=L(Q_\om(x),\widetilde Q(x), 0, \be)\in W$
be chosen such that  (16) holds. Let  $\{\Psi^\br{m}_{ni}(x)\}$  be
the solution of the main equation (25). The following assertion is
proved in [14].

\smallskip
{\bf Lemma 14. }{\it Consider the equations
$$
(I+A_0)y_0=f_0,\quad (I+A)y=f,
$$
in a Banach space $\mathfrak{B},$ where $A_0, A$ are linear bounded operators,
acting from $\mathfrak{B}$ to $\mathfrak{B},$ and $I$ is the identity operator.
Suppose that there exists the linear bounded operator $R_0:=(I+A_0)^{-1}.$
If $\|A-A_0\|\le (2\|R_0\|)^{-1},$  then there exists the linear bounded operator
$R=(I+A)^{-1},$ and $\|R\|\le2\|R_0\|,\;\;\;\;\|R-R_0\|\le2\|R_0\|^2\|A-A_0\|.$}

\smallskip
{\bf Lemma 15. }{\it The following relations hold
$$
\Psi^\br{m}_{ni}(x)\in C(\Om_\veps),\;\;|\Psi^\br{m}_{ni}(x)|\le C_\veps,  \eqno(37)
$$
$$
|\Psi^\br{m}_{ni}(x)-\widetilde\Psi^\br{m}_{ni}(x)|\le C_\veps\La\te_n,
\;\;\te_n=\left(\sum_{k=-\iy}^{+\iy}\frac{1}{(1+|\la_n^0-\la^0_k|)^2
(1+|\la^0_k|)^2}\right)^{1/2},\;\;x\in\Om_\veps,                           \eqno(38)
$$
$$
|(\Psi^\br{m}_{ni}(x))'|\le C_\veps(1+|\la_k^0|)\;\;x\in\Om_\veps,\;
|(\Psi^\br{m}_{ni}(x))'-(\widetilde\Psi^\br{m}_{ni}(x))'|
\le C_\veps\La,\;\; x\in\Om_\veps.                                         \eqno(39)
$$

Proof. } Using (26), we infer
$$
|\widetilde H_{ni,kj}(x)-\widetilde H_{ni,kj}(x_0)|\le|
\widetilde H'_{ni,kj}(\xi)|\,|x-x_0|\le
C_\veps|a_{k1}|\xi_k|x-x_0|,\;\;\;x,x_0,\xi\in\Om_\veps,
$$
hence
$$
\|\widetilde H(x)-\widetilde H(x_0)\|\le C_\veps\sum_{k=-\iy}^{+\iy}
|a_{k1}|\xi_k|x-x_0|\le C_\veps\La|x-x_0|.
$$
Choose $\de_0>0$ such that $\|\widetilde H(x)-\widetilde H(x_0)\|
\le (2\|\widetilde R(x_0)\|)^{-1}$ for $|x-x_0|\de_0$, where
$\widetilde R(x)=(I-\widetilde H(x))^{-1}.$ By lemma 14,
$\|\widetilde R(x)-\widetilde R(x_0)\|\le\
|\widetilde R(x_0)\|C_\veps|x-x_0|,\;\;x\in\Om_\veps.$
Therefore, the function $f(x)=\|\widetilde R(x)\|$ is continuous
and bounded on $\Om_\veps.$ Then
$$
\|\widetilde R(x)\|\le C_\veps,\;\;\;
\|\widetilde R(x)-\widetilde R(x_0)\|\le C_\veps|x-x_0|,\;\;x\in\Om_\veps,
$$
where the constant $C_\veps$ does not depend on $x, x_0$.
Since $\Psi^\br{m}(x)=\widetilde R(x)\widetilde\Psi^\br{m}(x),$
it follows that $\|\Psi^\br{m}(x)\|\le
\|\widetilde R(x)\|\,\|\widetilde\Psi^\br{m}(x)\|\le C_\veps.$
Thus, (37) is proved.

Taking (25)-(26) into account, we obtain
$$
|\Psi^\br{m}_{ni}(x)-\widetilde\Psi^\br{m}_{ni}(x)|\le
\sum_{k=-\iy}^{+\iy}\sum_{j=0}^1|\widetilde H_{ni,kj}(x)|\,|\Psi^\br{m}_{kj}(x)|
\le C_\veps\sum_{k=-\iy}^{+\iy}\frac{|a_{k1}|\xi_k}{1+|\la^0_k-\la^0_n|},
$$
and consequently,
$$
|\Psi^\br{m}_{ni}(x)-\widetilde\Psi^\br{m}_{ni}(x)|\le
C_\veps\te_n\sum_{k=-\iy}^{+\iy}\Big(|a_{k1}|\xi_k(1+|\la^0_k|)\Big)^2,
$$
i.e. (38) holds.
Estimates (39) are proved similarly. Lemma is proved.

\smallskip
Construct the functions $\vfi_{2,ni}(x)=\Big(\vfi_{12,ni}(x),\;
\vfi_{22,ni}(x)\Big)^T$ by the formula
$$
\vfi_{m2,n0}(x)=\Psi^\br{m}_{n0}(x)\xi_n+\Psi^\br{m}_{n1}(x),\;\;
\vfi_{m2,n1}(x)=\Psi^\br{m}_{n1}(x),\;\;m=1,2.                        \eqno(40)
$$
It follows from (40) and lemma 15 that
$$
|\vfi^{(m)}_{2,ni}(x)|\le C_\veps(1+|\la^0_n|)^m,\;\;\;
|\vfi^{(m)}_{2,n0}(x)-\vfi^{(m)}_{2,n1}(x)|
\le C_\veps\xi_n(1+|\la^0_n|)^m,\;\; m=0,1,                           \eqno(41)
$$
$$
|\widetilde\vfi_{2,ni}(x)-\vfi_{2,ni}(x)|\le C_\veps\La\te_n,
\;\;\; |\widetilde\vfi\,'_{2,ni}(x)-\vfi\,'_{2,ni}(x)|\le C_\veps\La.   \eqno(42)
$$

\smallskip
{\bf Lemma 16. }{\it The function $Q(x),$ constructed by (28),
is absolutely continuous on $[0,\pi].$

\smallskip
Proof. }
In view of  (41), the series in (29) converges uniformly on $\Om_\veps.$
According to lemma 15, the functions $\vfi_{2,ni}(x)$ are continuous,
and consequently, $\kp(x)$ is continuous on $\Om_\veps.$
One has
$$
\kp(x)=A_1(x)+A_2(x),\; A_1(x)=
\sum_{k=-\iy}^{+\iy}(a_{k0}-a_{k1})\widetilde\vfi_{2,k0}(x)\vfi^T_{2,k0}(x),
$$
$$
A_2(x)=\sum_{k=-\iy}^{+\iy}a_{k1}\Big(\widetilde\vfi_{2,k0}(x)\vfi^T_{2,k0}(x)
-\widetilde\vfi_{2,k1}(x)\vfi^T_{2,k1}(x)\Big).
$$
Taking (41) and (17) into account, we infer
$|(a_{k0}-a_{k1})\widetilde\vfi_{2,k0}(x)\vfi^T_{2,k0}(x)|
\le|a_{k1}|\xi_kC_\veps(1+|\la^0_k|).$
This yields that the series for $A_1(x)$ converges uniformly on
$\Om_\veps$ and $A'_1(x)\in L(0,\pi).$. Since
$$
(\widetilde\vfi_{2,k0}(x)\vfi^T_{2,k0}(x)-
\widetilde\vfi_{2,k1}(x)\vfi^T_{2,k1}(x))'=
(\widetilde\vfi\,'_{2,k0}(x)-\widetilde\vfi\,'_{2,k1}(x))\vfi^T_{2,k0}(x)
$$
$$
+\widetilde\vfi\,'_{2,k1}(x)(\vfi_{2,k0}(x)-\vfi_{2,k1}(x))^T+
 \widetilde\vfi_{2,k0}(x)(\vfi\,'_{2,k0}(x)-\vfi\,'_{2,k1}(x))^T
+(\widetilde\vfi_{2,k0}(x)-\widetilde\vfi_{2,k1}(x))
(\vfi^T_{2,k1}(x))',
$$
it follows from (41) and (17) that
$|a_{k1}(\widetilde\vfi_{2,k0}(x)\vfi^T_{2,k0}(x)-
\widetilde\vfi_{2,k1}(x)\vfi^T_{2,k1}(x))'|
\le C_\veps|a_{k1}|\xi_k(1+|\la^0_k|).$
This yields $A'_2(x)\in L(0,\pi).$ Thus, $\kp(x)$
is absolutely continuous on $[0,\pi]$ and $\kp'(x)\in L(0,\pi).$
Lemma is proved.

\smallskip
Let us now show that the given numbers $\{\la_k\}_{k=-\iy}^{+\iy}$
are eigenvalues of the constructed boundary value problem
$L(Q_\om(x),\;Q(x),\;0,\;\be)$.

\smallskip
{\bf Lemma 17. }{\it The following relations hold
$$
\ell\vfi_{2,kj}(x)=\la_{kj}\vfi_{2,kj}(x),
\;\;\ell\Phi_j(x,\la)=\la\Phi_j(x,\la),                                 \eqno(43)
$$
$$
\Phi_2(0,\la)=V_2(0),\;\;V^T_1(0)\Phi_1(0,\la)=1,\;\;
V^T_1(\be)\Phi_1(\pi,\la)=0,\;\;\De_{12}(\la_k)=0.                     \eqno(44)
$$

Proof.}
1) We construct $\Phi_j(x,\la)$ by (22). In view of (41)-(42),
the series in (22) converges uniformly in $\Om_\veps$.
Moreover, differentiating  (22) and taking (21) into account, we obtain
$$
\Phi'_j(x,\la)=\widetilde\Phi'_j(x,\la)+\sum_{k=-\iy}^{+\iy}
\Big(\widetilde D^\br{j}_{k0}(x,\la)a_{k0}\vfi\,'_{2,k0}(x)-
\widetilde D^\br{j}_{k1}(x,\la)a_{k1}\vfi\,'_{2,k1}(x)\Big)
$$
$$
-\sum_{k=-\iy}^{+\iy}\Big(\widetilde\Phi^T_j(x,\la)
\widetilde\vfi_{2,k0}(x)a_{k0}\vfi_{2,k0}(x)-
\widetilde\Phi^T_j(x,\la)\widetilde\vfi_{2,k1}(x)a_{k1}\vfi_{2,k1}(x)\Big),
$$
and consequently,
$$
(\Phi'_j(x,\la))^T=(\widetilde\Phi'_j(x,\la))^T-\widetilde\Phi^T_j(x,\la)\kp(x)
$$
$$
+\sum_{k=-\iy}^{+\iy}\Big(\widetilde D^\br{j}_{k0}(x,\la)a_{k0}
(\vfi\,'_{2,k0}(x))^T-\widetilde D^\br{j}_{k1}(x,\la)a_{k1}(\vfi\,'_{2,k1}(x))^T\Big).
$$
This yields
$$
\Big(B\Phi'_j(x,\la)+Q(x)\Phi_j(x,\la)\Big)^T=(\widetilde\Phi'_j(x,\la))^T
+\widetilde\Phi^T_j(x,\la)\Big(Q(x)+\kp(x)B\Big)+
$$
$$
\sum_{k=-\iy}^{+\iy}\Big(\widetilde
D^\br{j}_{k0}(x,\la)a_{k0}\Big(B\vfi\,'_{2,k0}(x)+Q(x)\vfi_{2,k0}(x)\Big)^T-\widetilde
D^\br{j}_{k1}(x,\la)a_{k1}\Big(B\vfi\,'_{2,k1}(x)+Q(x)\vfi_{2,k1}(x)\Big)^T\Big).
$$
Since $Q(x)=\widetilde Q(x)+B\kp(x)-\kp(x)B$,
$\widetilde\Phi^T_j(x,\la)B\widetilde\vfi_{2,ni}(x)
=\widetilde D^\br{j}_{ni}(x,\la)(\la-\la_{ni})$,  it follows that
$$
\widetilde\Phi^T_j(x,\la)B\kp(x)=\sum_{k=-\iy}^{+\iy}\Big(\widetilde
D^\br{j}_{k0}(x,\la)(\la-\la_{k0})a_{k0}\vfi^T_{2,k0}(x)-\widetilde
D^\br{j}_{k1}(x,\la)(\la-\la_{k1})a_{k1}\vfi^T_{2,k1}(x)\Big),
$$
hence
$$
\Big(B\Phi'_j(x,\la)+Q(x)\Phi_j(x,\la)\Big)^T=\Big(\widetilde\Phi'_j(x,\la)
+\widetilde Q(x)\widetilde\Phi_j(x,\la)\Big)^T
$$
$$
+\sum_{k=-\iy}^{+\iy}\Big(\widetilde D^\br{j}_{k0}(x,\la)a_{k0}
\Big(B\vfi\,'_{2,k0}(x)+Q(x)\vfi_{2,k0}(x)+(\la-\la_{k0})\vfi_{2,k0}(x))\Big)^T
$$
$$
-\widetilde D^\br{j}_{k1}(x,\la)a_{k1}\Big(B\vfi\,'_{2,k1}(x)+Q(x)\vfi_{2,k1}(x)
+(\la-\la_{k1})\vfi_{2,k1}(x)\Big)^T\Big).
$$
Taking (22) into account, we calculate
$$
\ell\Phi_j(x,\la)-\la\Phi_j(x,\la)=\sum_{k=-\iy}^{+\iy}
\Bigg(\widetilde D^\br{j}_{k0}(x,\la)a_{k0}\Big(\ell\vfi_{2,k0}(x)-
\la_{k0}\vfi_{2,k0}(x)\Big)
$$
$$
-\widetilde D^\br{j}_{k1}(x,\la)a_{k1}\Big(\ell\vfi_{2,k1}(x)
-\la_{k1}\vfi_{2,k1}(x)\Big)\Bigg).                                        \eqno(45)
$$
Consider (45) for $j=2$ and $\la=\la_{ni}$:
$$
z_{ni}(x)-\sum_{k=-\iy}^{+\iy}\Big(\widetilde P_{ni,k0}(x)z_{k0}(x)
-\widetilde P_{ni,k1}(x)z_{k1}(x)\Big)=0,
$$
where $z_{ni}(x)=\ell\vfi_{2,ni}(x)-\la_{ni}\vfi_{2,ni}(x),$ or
$$
Z^\br{m}_{ni}(x)-\sum_{k,j}\widetilde H_{ni,kj}(x)Z^\br{m}_{kj}(x)=0,     \eqno(46)
$$
where $Z^\br{m}_{n0}=\Big(z_{m,n0}(x)-z_{m,n1}(x)\Big)\chi_n,\;
Z^\br{m}_{n1}(x)=z_{m,n1}(x)$,
$z_{ni}(x)=(z_{1,ni}(x),\;z_{2,ni}(x))^T$, $m=1,2$.
taking (41) into account, we get  $|Z^\br{m}_{ni}(x)|\le
C_\veps(1+|\la^0_n|)$. Using (46) and (26), we infer
$$
|Z^\br{m}_{ni}(x)|\le C_\veps\sum_{k=-\iy}^{+\iy}
\frac{|a_{k1}|\xi_k(|\la^0_k|+1)}{1+|\la^0_n-\la^0_k|}\le C_\veps\La,
$$
and consequently, $\{Z^\br{m}_{ni}(x)\}\in\mathfrak{m}$. Equation
(46) has only trivial solution, i.e. $Z^\br{m}_{ni}(x)=0,$ hence
$\ell\vfi_{2,ni}(x)-\la_{ni}\vfi_{2,ni}(x)=0.$ The second relation (43)
follows now from (45).

2) Since $\widetilde D^\br{2}_{kj}(x,\la)=\di\frac{1}{\la-\la_{kj}}
\det(\widetilde\vfi_2(x,\la),\widetilde\vfi_{2,kj}(x)),$ it follows that
$$
\widetilde D^\br{2}_{kj}(0,\la)=\di\frac{1}{\la-\la_{kj}}\det(V_2(0),V_2(0))=0.
$$
Using (22), we find $\Phi_2(0,\la)=\widetilde\Phi_2(0,\la)=V_2(0).$
Furthermore, taking $j=1$, $x=0$ in (22) and multiplying by $V^T_1(0)$,
we calculate
$$
V^T_1(0)\Phi_1(0,\la)=V^T_1(0)\widetilde\Phi_1(0,\la)+\sum_{k=-\iy}^{+\iy}
\Big(\widetilde D^\br{1}_{k0}(0,\la)a_{k0}V^T_1(0)\vfi_{2,k0}(0)-
\widetilde D^\br{1}_{k1}(0,\la)a_{k1}V^T_1(0)\vfi_{2,k1}(0)\Big).
$$
Since $V^T_1(0)\vfi_{2,kj}(0)=V^T_1(0)V_2(0)=0,$ one gets
$V^T_1(0)\Phi_1(0,\la)=V^T_1(0)\widetilde\Phi_1(0,\la)=1$.
Thus, $\Phi_2(x,\la)=\vfi_2(x,\la)$ is a solution of  (1) with the
initial condition $\vfi_2(0,\la)=V_2(0).$ Then
$\De_{12}(\la)=V^T_1(\be)\vfi_2(\pi,\la).$
Let us show that $\De_{12}(\la_{n0})=0,$ i.e. $\{\la_n\}_{n=-\iy}^{+\iy}$
are eigenvalues of $L.$ For this purpose we take $j=2$, $x=\pi$
in (22) and multiply by $V^T_1(\be).$ This yields
$$
\De_{12}(\la)=\widetilde\De_{12}(\la)+\sum_{k=-\iy}^{+\iy}\Big(\widetilde
D^\br{2}_{k0}(\pi,\la)a_{k0}\De_{12}(\la_{k0})-\widetilde
D^\br{2}_{k1}(\pi,\la)a_{k1}\De_{12}(\la_{k1})\Big),
$$
and consequently,
$$
\De_{12}(\la_{ni})=\widetilde\De_{12}(\la_{ni})+\sum_{k=-\iy}^{+\iy}
\Big(\widetilde P_{ni,k0}(\pi)\De_{12}(\la_{k0})
-\widetilde P_{ni,k1}(\pi)\De_{12}(\la_{k1})\Big).                        \eqno(47)
$$
By definition $\De(\la)=V^T(\be)\vfi(\pi,\la)$, then
$\det\De(\la)\equiv1,$ or
$$
\De_{11}(\la)\De_{22}(\la)-\De_{12}(\la)\De_{21}(\la)\equiv1.             \eqno(48)
$$
Furthermore, $\br{\widetilde\vfi_2(\pi,\la),\,\widetilde\vfi_{2,kj}(\pi)}
=\widetilde\vfi^T_2(\pi,\la)B\widetilde\vfi_{2,kj}(\pi).$ Since
$V(\be)V^T(\be)=I,$ it follows that
$\bbr{\widetilde\vfi_2(\pi,\la),\,\widetilde\vfi_{2,kj}(\pi)}=
\Big(V^T(\be)\vfi_2(\pi,\la)\Big)^TB V^T(\be)\vfi_2(\pi,\la_{kj})$,
and consequently,
$$
\bbr{\widetilde\vfi_2(\pi,\la),\,\widetilde\vfi_{2,kj}(\pi)}=
\widetilde\De_{12}(\la)\widetilde\De_{22}(\la_{kj})-
\widetilde\De_{12}(\la_{kj})\widetilde\De_{22}(\la).
$$
Thus,
$$
\widetilde
D^\br{2}_{kj}(\pi,\la)=\frac{1}{\la-\la_{kj}}
\Big(\widetilde\De_{12}(\la)\widetilde\De_{22}(\la_{kj})-
\widetilde\De_{12}(\la_{kj})\widetilde\De_{22}(\la)\Big).              \eqno(49)
$$
From (49) for $n\neq k$, we find
$$
\widetilde P_{n1,k1}(\pi)=\di\frac{a_{k1}}{\la_{n1}-\la_{k1}}
\Big(\widetilde\De_{12}(\la_{n1})\widetilde\De_{22}(\la_{k1})-
\widetilde\De_{12}(\la_{k1})\widetilde\De_{22}(\la_{n1})\Big)=0.
$$
For $n=k,$ one has $\widetilde P_{n1,n1}(\pi)=a_{n1}
\dot{\widetilde\De}_{12}(\la_{n1})\widetilde\De_{22}(\la_{n1})$,
where $\dot{\widetilde\De}_{12}(\la):=\frac{d}{d\la}{\widetilde\De}_{12}(\la).$
Since $a_{n1}=\Res_{\la=\la_{n1}}\widetilde\fv(\la)=
-(\widetilde\De_{11}(\la_{n1}))(\dot{\widetilde\De}_{12}(\la_{n1}))^{-1}$,
it follows that $\widetilde P_{n1,n1}(\pi)=-\widetilde\De_{11}(\la_{n1})
\widetilde\De_{22}(\la_{n1})$. From (48) for $\la=\la_{n1}$ we infer
$\widetilde P_{n1,n1}(\pi)=-1.$ Thus,
$$
\widetilde P_{n1,k1}(\pi)=-\de_{nk},                                  \eqno(50)
$$
where $\de_{nk}$ is the Kronecker symbol.
From (49) for $\la_{n0}\neq\la_{k1}$ one has
$$
\widetilde P_{n0,k1}(\pi)=\di\frac{a_{k1}}{\la_{n0}-\la_{k1}}
\Big(\widetilde\De_{12}(\la_{n0})\widetilde\De_{22}(\la_{k1})-
\widetilde\De_{12}(\la_{k1})\widetilde\De_{22}(\la_{n0})\Big)=
a_{k1}\frac{\widetilde\De_{12}(\la_{n0})
\widetilde\De_{22}(\la_{k1})}{\la_{n0}-\la_{k1}}.
$$
By virtue of (48), $\widetilde\De_{22}(\la_{k1})=
(\widetilde\De_{11}(\la_{k1}))^{-1}$. Moreover,
$\widetilde P_{n0,k1}(\pi)=-1$ for $\la_{n0}=\la_{k1}$. Thus,
$$
\widetilde P_{n0,k1}(\pi)=-1\mbox{ for }\la_{n0}=\la_{k1},\;\;
\widetilde P_{n0,k1}(\pi)=-\di\frac{\widetilde\De_{12}(\la_{n0})}
{\dot{\widetilde\De}_{12}(\la_{k1})(\la_{n0}-\la_{k1})}
\mbox{ for }\la_{n0}\neq\la_{k1}.                                     \eqno(51)
$$

Consider the function $Z(\la)=(\De_{12}(\la)-\widetilde\De_{12}(\la))
(\widetilde\De_{12}(\la))^{-1}$. Since $Q_\om(x)=\widetilde Q_\om(x)$,
$\al=\widetilde\al=0$, $\be=\widetilde\be=0,$ it follows that
$\De_{12}(\la)-\widetilde\De_{12}(\la)=O\Big(e^{\pi|Im\la|}|\la|^{-\nu}\Big),$
hence $|Z(\la)|\le C_\de|\la|^{-\nu}$ for $\la\in\widetilde G_\de.$
In particular, this yields $\di\int_{|\xi|=\widetilde R_n}
\frac{Z(\xi)}{\xi-\la}d\xi\to0$ as $n\to\iy.$ Calculating the integral by
residue's theorem, we obtain for $n\to\iy,$
$$
\De_{12}(\la)=\widetilde\De_{12}(\la)+\sum_{k=-\iy}^{+\iy}
\frac{\widetilde\De_{12}(\la)}{(\la-\la_{k1})
\dot{\widetilde\De}_{12}(\la_{k1})} \De_{12}(\la_{k1}).
$$
Putting $\la=\la_{n0}$ and taking (51) into account, we infer
$$
\De_{12}(\la_{n0})=\widetilde\De_{12}(\la_{n0})-
\sum_{k=-\iy}^{+\iy}\widetilde P_{n0,k1}(\pi)\De_{12}(\la_{k1}).
$$
Together with (47) and (50) this yields
$$
\sum_{k=-\iy}^{+\iy}\widetilde P_{ni,k0}(\pi)\De_{12}(\la_{k0})=0,
$$
and consequently, $\De_{12}(\la_{n0})=0.$
Now we take $j=1$, $x=\pi$ in (22) and multiply by $V^T_1(\be).$ Then
$$
V^T_1(\be)\Phi_1(\pi,\la)=V^T_1(\be)\widetilde\Phi_1(\pi,\la)
+\sum_{k=-\iy}^{+\iy}\Big(\widetilde D^\br{1}_{k0}(\pi,\la)a_{k0}
\De_{12}(\la_0)-\widetilde D^\br{1}_{k1}(\pi,\la)a_{k1}\De_{12}(\la_1)\Big).
$$
Furthemore, $\widetilde D^\br{1}_{k1}(\pi,\la)=\di\frac{1}{\la-\la_{k1}}
(V^T(\be)\widetilde\Phi_1(\pi,\la))^T BV^T(\be)\widetilde\vfi_2(\pi,\la_{k1})$ or
$$
\widetilde D^\br{1}_{k1}(\pi,\la)=\frac{1}{\la-\la_{k1}}
\Big(V^T_1(\be)\widetilde\Phi_1(\pi,\la)\widetilde\De_{22}(\la_{k1})-
V^T_2(\be)\widetilde\Phi_1(\pi,\la)\widetilde\De_{12}(\la_{k1})\Big),
$$
hence, $\widetilde D^\br{1}_{k1}(\pi,\la)=0.$ Thus,
$V^T_1(\be)\Phi_1(\pi,\la)=V^T_1(\be)\widetilde\Phi_1(\pi,\la)=0$,
and all relations (44) are valid. Lemma 17 is proved.

\medskip
It remains to show that $a_k=\Res_{\la=\la_k}\fv(\la),$ where
$\fv(\la)=V^T_2(0)\Phi_1(0,\la).$ Taking in (22) $j=1$, $x=0$ and
multiplying by $V^T_2(0),$ we obtain
$$
\fv(\la)=\widetilde\fv(\la)+\sum_{k=-\iy}^{+\iy}
\Big(\widetilde D^\br{1}_{k0}(0,\la)a_{k0}V^T_2(0)\vfi_{2,k0}(0)-
\widetilde D^\br{1}_{k1}(0,\la)a_{k1}V^T_2(0)\vfi_{2,k1}(0)\Big).           \eqno(52)
$$
Using lemma 17, we calculate $V^T_2(0)\vfi_{2,kj}(0)=V^T_2(0)V_2(0)=1$.
Since
$$
\widetilde D^\br{1}_{kj}(0,\la)=\di\frac{1}{\la-\la_{kj}}
\widetilde\Phi^T_1(0,\la)B\widetilde\vfi_{2,kj}(0),\quad
\widetilde\vfi_{2,kj}(0)=V_2(0),
$$
it follows that
$\widetilde D^\br{1}_{kj}(0,\la)=\di\frac{1}{\la-\la_{kj}}$.
Thus, (52) takes the form
$$
\fv(\la)=\widetilde\fv(\la)+\sum_{k=-\iy}^{+\iy}\left(\frac{a_{k0}}{\la-\la_{k0}}
-\frac{a_{k1}}{\la-\la_{k1}}\right).
$$
With the help of Lemma 6, this yields $a_k=\Res_{\la=\la_k}\fv(\la).$
Theorem 3 is proved

\bigskip
{\bf Acknowledgment.} This work was supported by Grant 1.1436.2014K of the Russian
Ministry of Education and Science and by Grant 13-01-00134 of Russian Foundation for
Basic Research.

\begin{center}
{\bf REFERENCES}
\end{center}
\begin{enumerate}
\item[{[1]}] Meschanov V.P. and Feldstein A.L., Automatic Design of
     Directional Couplers, Moscow: Sviaz, 1980 (in Russian).
\item[{[2]}] Litvinenko O.N. and Soshnikov V.I., The Theory of
     Heterogenious Lines and their Applica\-tions in Radio Engineering,
     Moscow: Radio, 1964 (in Russian).
\item[{[3]}] Freiling G. and Yurko V.A., Reconstructing parameters
     of a medium from incomplete spectral information. Results in
     Mathematics 35 (1999), 228-249.
\item[{[4]}] Lapwood F.R. and Usami T., Free Oscilations of the Earth,
     Cambridge University Press, Cambridge, 1981.
\item[{[5]}] Wasow W. Linear turning point theory. Applied Mathematical
     Sciences, 54; Springer-Verlag, New York-Berlin, 1985.
\item[{[6]}] Gasymov M.G. Determination of Sturm-Liouville equation with
     a singular point from two spectra. Doklady Akad. Nauk SSSR 161 (1965),
     274-276; transl. in Sov. Math. Dokl. 6(1965), 396-399.
\item[{[7]}] Zhornitskaya L.A. and Serov V.S., Inverse eigenvalue
     problems  for  a singular  Sturm-Liouville  operator  on  (0,1).
     Inverse Problems 10 (1994), no.4, 975-987.
\item[{[8]}]  Yurko V.A., Inverse problem for differential equations
     with a singularity. Differen. Uravne\-niya 28 (1992), 1355-1362;
     English transl. in Differential Equations  28 (1992), 1100-1107.
\item[{[9]}] Yurko V.A., On higher-order differential operators with a
     singular point. Inverse Problems 9 (1993), 495-502.
\item[{[10]}] Yurko V.A., On higher-order differential operators with a
     regular  singularity.  Mat. Sb. 186 (1995),  no.6,  133-160;
     English  transl. in  Sbornik; Mathematics  186 (1995), no.6, 901-928.
\item[{[11]}] Yurko V.A., Integral transforms connected with
     differential operators having singularities inside the interval.
     Integral Transforms and Special Functions 5 (1997), no.3-4, 309-322.
\item[{[12]}] Gorbunov O.B.; Yurko V.A.; Shieh, Chung-Tsun. Spectral
     analysis of the Dirac system with a singularity in an interior
     point. arXiv:1410.2020v1 [math.SP], 17pp
\item[{[13]}] Yurko V A, Method of Spectral Mappings in the Inverse
     Problem Theory, Inverse and Ill-posed Problems Series, VSP, Utrecht, 2002.
\item[{[14]}] Freiling G.; Yurko V.A. Inverse Sturm-Liouville Problems and
     their Applications, Nova Science Publishers, Inc., Huntington, NY, 2001, 356 pp.
\item[{[15]}] Bellmann R. and Cooke K. Differential-difference Equations.
    Academic Press, New York, 1963.
\end{enumerate}

\begin{tabular}{ll}
Name:             &   Yurko, Viacheslav and Gorbunov, Oleg \\
Place of work:    &   Department of Mathematics, Saratov State University \\
{}                &   Astrakhanskaya 83, Saratov 410012, Russia \\
E-mail:           &   yurkova@info.sgu.ru  \\
\end{tabular}

\end{document}